%% file: ACM.tex
\documentclass[10pt,reqno]{amsart}
\usepackage[T1]{fontenc} 
\usepackage[utf8]{inputenc}               
\usepackage[english]{babel}                    
\usepackage{amsmath}
\usepackage{amsfonts}  
\usepackage{nicefrac}
\usepackage[pagebackref=true]{hyperref}
\usepackage{bm}
\usepackage[margin=2cm,includefoot,footskip=30pt]{geometry}
\usepackage{frontespizio}
\usepackage{eqparbox}
\usepackage{amssymb}
\usepackage{thmtools}
\usepackage{tikz-cd}
\usepackage{mathrsfs}
\usepackage{enumerate}
\usepackage{comment}
\usepackage{multicol}
\usepackage{longtable}
\usepackage{pdfpages}
\usepackage{array}
\usepackage[all]{xy}
\usepackage{listings}
\input{lst-Macaulay2}

\newtheorem{te}{Theorem}[section]
\newtheorem{pr}[te]{Proposition}
\newtheorem{lem}[te]{Lemma}
\newtheorem{co}[te]{Corollary}
\newtheorem{con}[te]{Conjecture}

\newtheorem{de}[te]{Definition}
\newtheorem{no}[te]{Notation}
\newtheorem{es}[te]{Example}
\newtheorem{Os}[te]{Remark}
\newtheorem{Qu}[te]{Question}

\newcommand{\N}{\mathbb{N}}
\newcommand{\Z}{\mathbb{Z}}
\newcommand{\Q}{\mathbb{Q}}

\newcommand{\C}{\mathbb{C}}

\newcommand{\p}{\mathbb{P}}
\newcommand{\os}{\mathcal{O}}
\newcommand{\osn}{\os_{\p^n}}
\newcommand{\e}{\mathcal{E}}
\newcommand{\f}{\mathcal{F}}
\newcommand{\g}{\mathcal{G}}
\newcommand{\id}{\mathcal{I}}

\renewcommand{\c}{\mathscr{C}}

\newcommand{\ck}{\mathcal{K}}
\newcommand{\cl}{\mathcal{L}}

\newcommand{\n}{\mathcal{N}}

\newcommand{\cu}{\mathcal{U}}

\newcommand{\gel}[1]{\underline{\textbf{#1)}}}
\newcommand{\du}{^\vee}
\newcommand{\dd}{^{\vee\vee}}

\newcommand{\mt}{\mapsto}

\newcommand{\inv}{^{-1}}

\newcommand{\gk}{\mathbf{k}}
\newcommand{\gi}{\mathbf{I}}

\newcommand{\sesl}[4]{\begin{equation}\label{#4} 0\to #1\to #2\to #3 \to 0\end{equation}}
\newcommand{\ses}[3]{\begin{equation*} 0\to #1\to #2\to #3 \to 0\end{equation*}}

\newcommand{\on}[1]{\os_{\mathbb{P}^{#1}}}

\newcommand{\ext}{\mathscr{E}xt}
\newcommand{\Hom}{\mathscr{H}om}

\newcommand{\De}{\begin{de}}
\newcommand{\Ne}{\end{de}}
\newcommand{\Te}{\begin{te}}
\newcommand{\Ma}{\end{te}}
\newcommand{\Prop}{\begin{pr}}
\newcommand{\One}{\end{pr}}
\newcommand{\Es}{\begin{es}}
\newcommand{\Io}{\end{es}}
\newcommand{\Le}{\begin{lem}}
\newcommand{\ma}{\end{lem}}
\newcommand{\Co}{\begin{co}}
\newcommand{\io}{\end{co}}
\newcommand{\Con}{\begin{con}}
\newcommand{\ura}{\end{con}}
\newcommand{\Oss}{\begin{Os}}
\newcommand{\one}{\end{Os}}
\newcommand{\Do}{\begin{Qu}}
\newcommand{\da}{\end{Qu}}
\newcommand{\No}{\begin{no}}
\newcommand{\ione}{\end{no}}

\newcommand{\dia}{ \begin{center}\begin{tikzcd}}
\newcommand{\mma}{\end{tikzcd}\end{center}}

\hypersetup{colorlinks,
    linkcolor={green!80!blue},
    citecolor={blue!90!black},
    urlcolor={blue!70!black}
}

\address{\parbox{0.9\textwidth}{Dipartimento di Matematica Università degli Studi di Roma Tor Vergata \\[1pt]
Via della Ricerca Scientifica, 1 - 00133, Roma
\vspace{1mm}}}
\email{{vacca@mat.uniroma2.it}}
\title{Rank $2$ aCM and Ulrich bundles on Fano and Calabi--Yau double coverings of $\p^3$}

\author{Roberto Vacca}

\begin{document}

\begin{abstract}
We prove existence of aCM and Ulrich sheaves respect to ample and globally generated polarisations on a class of special finite coverings $f:X\to\p^n$, which in particular contains cyclic ones.
In the case of rank $2$ on double coverings, we have a precise description of the zero loci of such sheaves which allows us to study their geometry and classify all possible such bundles in the case $X$ is regular. 
We show that on a general double covering of $\p^3$ branched along a divisor of degree $4,6,8$ all the above sheaves exist and, when stable, we compute the dimension of their component in the moduli spaces.
\end{abstract}

\maketitle

\tableofcontents

\input{intro}
\input{generali}
\input{double}

\newpage

\bibliographystyle{alphaurl}
\bibliography{bibliografia}

\end{document}

%% file: lst-Macaulay2.tex
\usepackage{listings}
\usepackage{xcolor}
\definecolor{maccolor}{rgb}{0.3,0.3,0.8}
\lstdefinelanguage{Macaulay2}
{
basicstyle={\ttfamily},
keywordstyle={\color{maccolor!80!black}},
commentstyle={\color{gray}},
stringstyle={\color{red!40!black}},
rulecolor=\color{maccolor},
basewidth={1.2ex}, 
sensitive=false,
morecomment=[l]{--},
morecomment=[s]{-*}{*-},
morestring=[b]",
escapechar={`},
escapebegin={\rmfamily},
morekeywords={about,abs,AbstractToricVarieties,accumulate,Acknowledgement,acos,acosh,acot,addCancelTask,addDependencyTask,addEndFunction,addHook,AdditionalPaths,addStartFunction,addStartTask,Adjacent,adjoint,AdjointIdeal,AffineVariety,AfterEval,AfterNoPrint,AfterPrint,agm,AInfinity,alarm,AlgebraicSplines,Algorithm,Alignment,all,AllCodimensions,allowableThreads,ambient,analyticSpread,Analyzer,AnalyzeSheafOnP1,ancestor,ancestors,ANCHOR,and,andP,AngleBarList,ann,annihilator,antipode,any,append,applicationDirectory,applicationDirectorySuffix,apply,applyKeys,applyPairs,applyTable,applyValues,apropos,argument,Array,arXiv,Ascending,ascii,asin,asinh,ass,assert,associatedGradedRing,associatedPrimes,AssociativeAlgebras,AssociativeExpression,atan,atan2,atEndOfFile,Authors,autoload,AuxiliaryFiles,backtrace,Bag,Bareiss,baseFilename,BaseFunction,baseName,baseRing,baseRings,BaseRow,BasicList,basis,BasisElementLimit,Bayer,BeforePrint,beginDocumentation,BeginningMacaulay2,Benchmark,benchmark,Bertini,BesselJ,BesselY,betti,BettiCharacters,BettiTally,between,BGG,BIBasis,Binary,BinaryOperation,Binomial,binomial,BinomialEdgeIdeals,Binomials,BKZ,BlockMatrix,BLOCKQUOTE,BODY,Body,BoijSoederberg,BOLD,Book3264Examples,Boolean,BooleanGB,borel,Boxes,BR,break,Browse,Bruns,cache,CacheExampleOutput,CacheFunction,CacheTable,cacheValue,CallLimit,cancelTask,capture,catch,Caveat,CC,CDATA,ceiling,Center,centerString,Certification,ChainComplex,chainComplex,ChainComplexExtras,ChainComplexMap,ChainComplexOperations,ChangeMatrix,char,CharacteristicClasses,characters,charAnalyzer,check,CheckDocumentation,chi,Chordal,class,Classic,clean,clearAll,clearEcho,clearOutput,close,closeIn,closeOut,ClosestFit,CODE,code,codim,CodimensionLimit,coefficient,CoefficientRing,coefficientRing,coefficients,Cofactor,CohenEngine,CohenTopLevel,CoherentSheaf,CohomCalg,cohomology,coimage,CoincidentRootLoci,coker,cokernel,collectGarbage,columnAdd,columnate,columnMult,columnPermute,columnRankProfile,columnSwap,combine,Command,commandInterpreter,commandLine,COMMENT,commonest,commonRing,comodule,CompactMatrix,compactMatrixForm,CompiledFunction,CompiledFunctionBody,CompiledFunctionClosure,Complement,complement,complete,CompleteIntersection,CompleteIntersectionResolutions,Complexes,ComplexField,components,compose,compositions,compress,concatenate,conductor,ConductorElement,cone,Configuration,ConformalBlocks,conjugate,connectionCount,Consequences,Constant,Constants,constParser,content,continue,contract,Contributors,ConvexInterface,conwayPolynomial,ConwayPolynomials,copy,copyDirectory,copyFile,copyright,Core,CorrespondenceScrolls,cos,cosh,cot,CotangentSchubert,cotangentSheaf,coth,cover,coverMap,cpuTime,createTask,Cremona,csc,csch,current,currentColumnNumber,currentDirectory,currentFileDirectory,currentFileName,currentLayout,currentLineNumber,currentPackage,currentString,currentTime,Cyclotomic,Database,Date,DD,dd,deadParser,debug,debugError,DebuggingMode,debuggingMode,debugLevel,DecomposableSparseSystems,Decompose,decompose,deepSplice,Default,default,defaultPrecision,Degree,degree,degreeLength,DegreeLift,DegreeLimit,DegreeMap,DegreeOrder,DegreeRank,Degrees,degrees,degreesMonoid,degreesRing,delete,demark,denominator,Dense,Density,Depth,depth,Descending,Descent,Describe,describe,Description,det,determinant,DeterminantalRepresentations,DGAlgebras,diagonalMatrix,diameter,Dictionary,dictionary,dictionaryPath,diff,DiffAlg,difference,dim,directSum,disassemble,discriminant,dismiss,Dispatch,distinguished,DIV,Divide,divideByVariable,DivideConquer,DividedPowers,Divisor,DL,Dmodules,do,doc,docExample,docTemplate,document,DocumentTag,Down,drop,DT,dual,eagonNorthcott,EagonResolution,echoOff,echoOn,EdgeIdeals,edit,EigenSolver,eigenvalues,eigenvectors,eint,EisenbudHunekeVasconcelos,elapsedTime,elapsedTiming,elements,Eliminate,eliminate,Elimination,EliminationMatrices,EllipticCurves,EllipticIntegrals,else,EM,Email,End,end,endl,endPackage,Engine,engineDebugLevel,EngineRing,EngineTests,entries,EnumerationCurves,environment,Equation,EquivariantGB,erase,erf,erfc,error,errorDepth,euler,EulerConstant,eulers,even,EXAMPLE,ExampleFiles,ExampleItem,examples,ExampleSystems,Exclude,exec,exit,exp,expectedReesIdeal,expm1,exponents,export,exportFrom,exportMutable,Expression,expression,Ext,extend,ExteriorIdeals,ExteriorModules,exteriorPower,Factor,factor,false,Fano,FastMinors,FastNonminimal,FGLM,File,fileDictionaries,fileExecutable,fileExists,fileExitHooks,fileLength,fileMode,FileName,FilePosition,fileReadable,fileTime,fileWritable,fillMatrix,findFiles,findHeft,FindOne,findProgram,findSynonyms,FiniteFittingIdeals,First,first,firstkey,FirstPackage,fittingIdeal,flagLookup,FlatMonoid,flatten,flattenRing,Flexible,flip,floor,flush,fold,FollowLinks,for,forceGB,fork,FormalGroupLaws,Format,format,formation,FourierMotzkin,FourTiTwo,fpLLL,frac,fraction,FractionField,frames,FrobeniusThresholds,from,fromDividedPowers,fromDual,Function,FunctionApplication,FunctionBody,functionBody,FunctionClosure,FunctionFieldDesingularization,fusePairs,futureParser,GaloisField,Gamma,gb,GBDegrees,gbRemove,gbSnapshot,gbTrace,gcd,gcdCoefficients,gcdLLL,GCstats,genera,GeneralOrderedMonoid,GenerateAssertions,generateAssertions,generator,generators,Generic,GenericInitialIdeal,genericMatrix,genericSkewMatrix,genericSymmetricMatrix,gens,genus,get,getc,getChangeMatrix,getenv,getGlobalSymbol,getNetFile,getNonUnit,getPrimeWithRootOfUnity,getSymbol,getWWW,GF,gfanInterface,Givens,GKMVarieties,GLex,Global,global,globalAssign,globalAssignFunction,GlobalAssignHook,globalAssignment,globalAssignmentHooks,GlobalDictionary,GlobalHookStore,globalReleaseFunction,GlobalReleaseHook,Gorenstein,GradedLieAlgebras,GradedModule,gradedModule,GradedModuleMap,gradedModuleMap,gramm,GraphicalModels,GraphicalModelsMLE,Graphics,graphIdeal,graphRing,Graphs,Grassmannian,GRevLex,GroebnerBasis,groebnerBasis,GroebnerBasisOptions,GroebnerStrata,GroebnerWalk,groupID,GroupLex,GroupRevLex,GTZ,Hadamard,handleInterrupts,HardDegreeLimit,hash,HashTable,hashTable,HEAD,HEADER1,HEADER2,HEADER3,HEADER4,HEADER5,HEADER6,HeaderType,Heading,Headline,Heft,heft,Height,height,help,Hermite,hermite,Hermitian,HH,hh,HigherCIOperators,HighestWeights,Hilbert,hilbertFunction,hilbertPolynomial,hilbertSeries,HodgeIntegrals,hold,Holder,Hom,homeDirectory,HomePage,Homogeneous,Homogeneous2,homogenize,homology,homomorphism,HomotopyLieAlgebra,hooks,horizontalJoin,HorizontalSpace,HR,HREF,HTML,html,httpHeaders,Hybrid,HyperplaneArrangements,Hypertext,hypertext,HypertextContainer,HypertextParagraph,icFracP,icFractions,icMap,icPIdeal,id,Ideal,ideal,idealizer,identity,if,IgnoreExampleErrors,ii,image,imaginaryPart,IMG,ImmutableType,importFrom,in,incomparable,Increment,independentSets,indeterminate,IndeterminateNumber,Index,index,indexComponents,IndexedVariable,IndexedVariableTable,indices,inducedMap,inducesWellDefinedMap,InexactField,InexactFieldFamily,InexactNumber,InfiniteNumber,infinity,info,InfoDirSection,infoHelp,Inhomogeneous,input,Inputs,insert,installAssignmentMethod,installedPackages,installHilbertFunction,installMethod,installMinprimes,installPackage,InstallPrefix,instance,instances,IntegralClosure,integralClosure,integrate,IntermediateMarkUpType,interpreterDepth,intersect,intersectInP,Intersection,intersection,interval,InvariantRing,inverse,InverseMethod,inversePermutation,Inverses,inverseSystem,InverseSystems,Invertible,InvolutiveBases,irreducibleCharacteristicSeries,irreducibleDecomposition,isAffineRing,isANumber,isBorel,isCanceled,isCommutative,isConstant,isDirectory,isDirectSum,isEmpty,isField,isFinite,isFinitePrimeField,isFreeModule,isGlobalSymbol,isHomogeneous,isIdeal,isInfinite,isInjective,isInputFile,isIsomorphism,isLinearType,isListener,isLLL,isMember,isModule,isMonomialIdeal,isNormal,isOpen,isOutputFile,isPolynomialRing,isPrimary,isPrime,isPrimitive,isPseudoprime,isQuotientModule,isQuotientOf,isQuotientRing,isReady,isReal,isReduction,isRegularFile,isRing,isSkewCommutative,isSorted,isSquareFree,isStandardGradedPolynomialRing,isSubmodule,isSubquotient,isSubset,isSupportedInZeroLocus,isSurjective,isTable,isUnit,isWellDefined,isWeylAlgebra,ITALIC,Iterate,Jacobian,jacobian,jacobianDual,Jets,Join,join,Jupyter,K3Carpets,K3Surfaces,Keep,KeepFiles,KeepZeroes,ker,kernel,kernelLLL,kernelOfLocalization,Key,keys,Keyword,Keywords,kill,koszul,Kronecker,KustinMiller,LABEL,last,lastMatch,LATER,LatticePolytopes,Layout,lcm,leadCoefficient,leadComponent,leadMonomial,leadTerm,Left,left,length,LengthLimit,letterParser,Lex,LexIdeals,LI,Licenses,LieTypes,lift,liftable,Limit,limitFiles,limitProcesses,Linear,LinearAlgebra,LinearTruncations,lineNumber,lines,LINK,linkFile,List,list,listForm,listLocalSymbols,listSymbols,listUserSymbols,LITERAL,LLL,LLLBases,lngamma,load,loadDepth,LoadDocumentation,loadedFiles,loadedPackages,loadPackage,Local,local,localDictionaries,LocalDictionary,localize,LocalRings,locate,log,log1p,LongPolynomial,lookup,lookupCount,LowerBound,LUdecomposition,M0nbar,M2CODE,Macaulay2Doc,makeDirectory,MakeDocumentation,makeDocumentTag,MakeHTML,MakeInfo,MakeLinks,makePackageIndex,MakePDF,makeS2,Manipulator,map,MapExpression,MapleInterface,markedGB,Markov,MarkUpType,match,mathML,Matrix,matrix,MatrixExpression,Matroids,max,maxAllowableThreads,maxExponent,MaximalRank,maxPosition,MaxReductionCount,MCMApproximations,member,memoize,memoizeClear,memoizeValues,MENU,merge,mergePairs,META,method,MethodFunction,MethodFunctionBinary,MethodFunctionSingle,MethodFunctionWithOptions,methodOptions,methods,midpoint,min,minExponent,mingens,mingle,minimalBetti,MinimalGenerators,MinimalMatrix,minimalPresentation,minimalPresentationMap,minimalPresentationMapInv,MinimalPrimes,minimalPrimes,minimalReduction,Minimize,minimizeFilename,MinimumVersion,minors,minPosition,minPres,minprimes,Minus,minus,Miura,MixedMultiplicity,mkdir,mod,Module,module,ModuleDeformations,modulo,MonodromySolver,Monoid,monoid,MonoidElement,Monomial,MonomialAlgebras,monomialCurveIdeal,MonomialIdeal,monomialIdeal,MonomialIntegerPrograms,MonomialOrbits,MonomialOrder,Monomials,monomials,MonomialSize,monomialSubideal,moveFile,multidegree,multidoc,multigraded,MultigradedBettiTally,MultiGradedRationalMap,multiplicity,MultiplicitySequence,MultiplierIdeals,MultiplierIdealsDim2,MultiprojectiveVarieties,mutable,MutableHashTable,mutableIdentity,MutableList,MutableMatrix,mutableMatrix,NAGtypes,Name,nanosleep,Nauty,NautyGraphs,NCAlgebra,NCLex,needs,needsPackage,Net,net,NetFile,netList,new,newClass,newCoordinateSystem,NewFromMethod,newline,NewMethod,newNetFile,NewOfFromMethod,NewOfMethod,newPackage,newRing,nextkey,nextPrime,nil,NNParser,NoetherianOperators,NoetherNormalization,NonAssociativeProduct,NonminimalComplexes,nonspaceAnalyzer,NoPrint,norm,normalCone,Normaliz,NormalToricVarieties,not,Nothing,notify,notImplemented,NTL,null,nullaryMethods,nullhomotopy,nullParser,nullSpace,Number,number,NumberedVerticalList,numcols,numColumns,numerator,numeric,NumericalAlgebraicGeometry,NumericalCertification,NumericalImplicitization,NumericalLinearAlgebra,NumericalSchubertCalculus,numericInterval,NumericSolutions,numgens,numRows,numrows,odd,oeis,of,ofClass,OL,OldPolyhedra,OldToricVectorBundles,on,OneExpression,OnlineLookup,OO,oo,ooo,oooo,openDatabase,openDatabaseOut,openFiles,openIn,openInOut,openListener,OpenMath,openOut,openOutAppend,operatorAttributes,Option,OptionalComponentsPresent,optionalSignParser,Options,options,OptionTable,optP,or,Order,order,OrderedMonoid,orP,OutputDictionary,Outputs,override,pack,Package,package,PackageCitations,PackageDictionary,PackageExports,PackageImports,PackageTemplate,packageTemplate,pad,pager,PairLimit,pairs,PairsRemaining,PARA,Parametrization,parent,Parenthesize,Parser,Parsing,part,Partition,partition,partitions,parts,path,pdim,peek,PencilsOfQuadrics,Permanents,permanents,permutations,pfaffians,PHCpack,PhylogeneticTrees,pi,PieriMaps,pivots,PlaneCurveSingularities,plus,poincare,poincareN,Points,polarize,poly,Polyhedra,Polymake,PolynomialRing,Posets,Position,position,positions,PositivityToricBundles,POSIX,Postfix,Power,power,powermod,PRE,Precision,precision,Prefix,prefixDirectory,prefixPath,preimage,prepend,presentation,pretty,primaryComponent,PrimaryDecomposition,primaryDecomposition,PrimaryTag,PrimitiveElement,Print,print,printerr,printingAccuracy,printingLeadLimit,printingPrecision,printingSeparator,printingTimeLimit,printingTrailLimit,printString,printWidth,processID,Product,product,ProductOrder,profile,profileSummary,Program,programPaths,ProgramRun,Proj,Projective,ProjectiveHilbertPolynomial,projectiveHilbertPolynomial,ProjectiveVariety,promote,protect,Prune,prune,PruneComplex,pruningMap,Pseudocode,pseudocode,pseudoRemainder,Pullback,PushForward,pushForward,Python,QQ,QQParser,QRDecomposition,QthPower,Quasidegrees,QuaternaryQuartics,QuillenSuslin,quit,Quotient,quotient,quotientRemainder,QuotientRing,Radical,radical,RadicalCodim1,radicalContainment,RaiseError,random,RandomCanonicalCurves,RandomComplexes,RandomCurves,RandomCurvesOverVerySmallFiniteFields,RandomGenus14Curves,RandomIdeals,randomKRationalPoint,RandomMonomialIdeals,randomMutableMatrix,RandomObjects,RandomPlaneCurves,RandomPoints,RandomSpaceCurves,Range,rank,RationalMaps,RationalPoints,RationalPoints2,ReactionNetworks,read,readDirectory,readlink,readPackage,RealField,RealFP,realPart,realpath,RealQP,RealQP1,RealRoots,RealRR,RealXD,recursionDepth,recursionLimit,Reduce,reducedRowEchelonForm,reduceHilbert,reductionNumber,ReesAlgebra,reesAlgebra,reesAlgebraIdeal,reesIdeal,References,ReflexivePolytopesDB,regex,regexQuote,registerFinalizer,regSeqInIdeal,Regularity,regularity,relations,RelativeCanonicalResolution,relativizeFilename,Reload,remainder,RemakeAllDocumentation,remove,removeDirectory,removeFile,removeLowestDimension,reorganize,replace,RerunExamples,res,reshape,ResidualIntersections,ResLengthThree,Resolution,resolution,ResolutionsOfStanleyReisnerRings,restart,Result,resultant,Resultants,return,returnCode,Reverse,reverse,RevLex,Right,right,Ring,ring,RingElement,RingFamily,ringFromFractions,RingMap,rootPath,roots,rootURI,rotate,round,rowAdd,RowExpression,rowMult,rowPermute,rowRankProfile,rowSwap,RR,RRi,rsort,run,RunDirectory,RunExamples,RunExternalM2,runHooks,runLengthEncode,runProgram,same,saturate,Saturation,scan,scanKeys,scanLines,scanPairs,scanValues,schedule,schreyerOrder,Schubert,Schubert2,SchurComplexes,SchurFunctors,SchurRings,SCRIPT,scriptCommandLine,ScriptedFunctor,SCSCP,searchPath,sec,sech,SectionRing,SeeAlso,seeParsing,SegreClasses,select,selectInSubring,selectVariables,SelfInitializingType,SemidefiniteProgramming,Seminormalization,separate,SeparateExec,separateRegexp,Sequence,sequence,Serialization,serialNumber,Set,set,setEcho,setGroupID,setIOExclusive,setIOSynchronized,setIOUnSynchronized,setRandomSeed,setup,setupEmacs,sheaf,SheafExpression,sheafExt,sheafHom,SheafOfRings,shield,ShimoyamaYokoyama,short,show,showClassStructure,showHtml,showStructure,showTex,showUserStructure,SimpleDoc,simpleDocFrob,SimplicialComplexes,SimplicialDecomposability,SimplicialPosets,SimplifyFractions,sin,singularLocus,sinh,size,size2,SizeLimit,SkewCommutative,SlackIdeals,sleep,SLnEquivariantMatrices,SLPexpressions,SMALL,smithNormalForm,solve,someTerms,Sort,sort,sortColumns,SortStrategy,source,SourceCode,SourceRing,SPACE,SpaceCurves,SPAN,span,SparseMonomialVectorExpression,SparseResultants,SparseVectorExpression,Spec,SpechtModule,SpecialFanoFourfolds,specialFiber,specialFiberIdeal,SpectralSequences,splice,splitWWW,sqrt,SRdeformations,stack,stacksProject,Standard,standardForm,standardPairs,StartWithOneMinor,stashValue,StatePolytope,StatGraphs,status,stderr,stdio,step,StopBeforeComputation,stopIfError,StopWithMinimalGenerators,Strategy,String,STRONG,StronglyStableIdeals,STYLE,Style,style,SUB,sub,SubalgebraBases,sublists,submatrix,submatrixByDegrees,Subnodes,subquotient,SubringLimit,Subscript,subscript,SUBSECTION,subsets,substitute,substring,subtable,Sugarless,Sum,sum,SumOfTwists,SumsOfSquares,SUP,super,SuperLinearAlgebra,Superscript,superscript,support,SVD,SVDComplexes,switch,SwitchingFields,sylvesterMatrix,Symbol,symbol,SymbolBody,symbolBody,SymbolicPowers,symlinkDirectory,symlinkFile,symmetricAlgebra,symmetricAlgebraIdeal,symmetricKernel,SymmetricPolynomials,symmetricPower,synonym,SYNOPSIS,syz,Syzygies,SyzygyLimit,SyzygyMatrix,SyzygyRows,syzygyScheme,TABLE,Table,table,take,Tally,tally,tan,TangentCone,tangentCone,tangentSheaf,tanh,target,Task,taskResult,TateOnProducts,TD,temporaryFileName,tensor,tensorAssociativity,TensorComplexes,terminalParser,terms,TEST,Test,testExample,testHunekeQuestion,TestIdeals,TestInput,tests,TEX,tex,TeXmacs,texMath,Text,TH,then,Thing,ThinSincereQuivers,ThreadedGB,threadVariable,Threshold,throw,Time,time,times,timing,TITLE,TO,to,TO2,toAbsolutePath,toCC,toDividedPowers,toDual,toExternalString,toField,TOH,toList,toLower,top,top,topCoefficients,Topcom,topComponents,topLevelMode,Tor,TorAlgebra,Toric,ToricInvariants,ToricTopology,ToricVectorBundles,toRR,toRRi,toSequence,toString,TotalPairs,toUpper,TR,trace,transpose,TriangularSets,Tries,Trim,trim,Triplets,Tropical,true,Truncate,truncate,truncateOutput,Truncations,try,TSpreadIdeals,TT,tutorial,Type,TypicalValue,typicalValues,UL,ultimate,unbag,uncurry,Undo,undocumented,uniform,uninstallAllPackages,uninstallPackage,Unique,unique,Units,Unmixed,unsequence,unstack,Up,UpdateOnly,UpperTriangular,URL,urlEncode,Usage,use,UseCachedExampleOutput,UseHilbertFunction,UserMode,userSymbols,UseSyzygies,utf8,utf8check,validate,value,values,Variable,VariableBaseName,Variables,Variety,variety,vars,Vasconcelos,Vector,vector,VectorExpression,VectorFields,VectorGraphics,Verbose,Verbosity,Verify,VersalDeformations,versalEmbedding,Version,version,VerticalList,VerticalSpace,viewHelp,VirtualResolutions,VirtualTally,VisibleList,Visualize,wait,WebApp,wedgeProduct,weightRange,Weights,WeylAlgebra,WeylGroups,when,whichGm,while,width,wikipedia,Wrap,wrap,WrapperType,XML,xor,youngest,zero,ZeroExpression,zeta,ZZ,ZZParser}
}
\lstalias{Macaulay2output}{Macaulay2}

%% file: intro.tex
\section{Introduction}
The aim of this paper is to study a special class of vector bundles, characterised by some cohomological vanishings.
Given a polarised variety $(X,H)$ of dimension $n$ over a field $\gk$, we define a coherent sheaf $\e$ to be arithmetically Cohen-Macaulay (aCM) if $h^j(\e(i)))=0$ for all $i\in \Z$ and $0<j<n$ and $depth(\e_x)=dim(\os_{X,x})$ for all $x\in X$.
The first fundamental result on aCM vector bundles due to Horrocks, later strengthened by Evans--Griffith, is the following.

\medskip
\noindent\textbf{Theorem} 
(\cite[Prop. 9.5]{Horrocks}, \cite[Thm. 2.4]{EvansGriffith})
\textit{A vector bundle $\e$ on $\p^n$ is aCM if and only if it splits as a direct sum of line bundles if and only if $h^j(\e(i)))=0$ for all $i\in \Z$ and $0<j<min\{rk(\e),n\}$.}
\medskip

aCM sheaves have been studied in the eighties first from the point of view of commutative algebra, since they correspond to maximal Cohen-Macaulay modules.
Remarkably, Eisenbud showed that aCM bundles (modules) on divisors identify matrix factorisations of their equation, see \cite{Eis-resolutions}, and later this result was extended to sheaves by Beauville \cite{Beadeterminantal}.
Matrix factorisations in turn are connected to singularity categories of hypersurfaces and mirror symmetry by Orlov \cite{Orlovmatrix}.
In the last decades, the construction and classification of aCM bundles on a variety $X$ have garnered significant attention, particularly in connection with the study of the bounded derived category of $X$, see \cite{Ottaviani_split} for a survey. 
The presence and the structure of the families of aCM bundles can be seen as a measure of complexity of our variety $X$.
From this point of view, we can distinguish $3$ cases: \textit{finite, tame, wild} depending on whether the parameter spaces for such sheaves are made of a finite set of points, a countable union of points and curves or varieties whose dimension becomes arbitrary high when the rank increases.
It has been proven by Faenzi and Pons-Llopis \cite{FaenziPons_acm} that a reduced aCM subscheme\footnote{$X\subset\p^N$ is aCM if and only if $X$ is projectively normal and $\os_X$ is an aCM sheaf} $X\subset \p^N$ is of wild type except for a very specific list of exceptions which are of finite or tame type.

Further, we can ask for the existence of even more special sheaves: Ulrich sheaves.
In this case, to the aCM condition we also add the vanishing of $h^{n}(\e(-n))$.
The corresponding modules are equivalently characterised as maximal Cohen-Macaulay modules with the maximum possible number of generators, see \cite{Brenna_Herzog_Ulrich}.
The existence of those sheaves implies that the structure of the cone of cohomology tables of sheaves on $X$ is the same as the one of $\p^n$, see \cite{Boij-Sod} and \cite{EisSch_BStheory}, and if $H$ is very ample, they detect the possibility of writing the Chow form of $X$ as a determinant of a matrix whose entries are linear forms, see \cite{EisSch}.
The main questions in the field have been asked by Eisenbud and Schreyer in the just cited work:

\smallskip
\noindent\textbf{Question}(Eisenbud--Schreyer). 
\textit{Does any projective variety $X\subset \p^N$ admit an Ulrich sheaf?
If yes, what is the minimal rank of such a sheaf?}
\smallskip

There is a great deal of literature on these topics, see, for example, the survey \cite{Bea}, the book \cite{CMRPL}.
The simplest possible case besides projective space are quadrics, on which indecomposable aCM bundles have been classified in \cite{Knorrer}: they are line bundles and spinor bundles; the latter ones are actually Ulrich.
Similarly, the case of complete intersections in projective spaces has been faced using commutative algebra methods, see \cite{Backelin_Herzog_Sanders} and \cite{BHU}.
Moreover, curves always admits Ulrich sheaves \cite[Cor. 4.5]{EisSch} and smooth complex surfaces always have rank $2$ Ulrich bundles, given that we choose a sufficiently large multiple of our polarisation \cite[Thm. 4.3]{CoskunHuizenga}.

\medskip

In this work we consider a variety $X$ with an ample and globally generated polarisation $H$.
Although most of the previous cited work assumed the polarisation to be very ample, almost all the properties of aCM and Ulrich sheaves are encoded in cohomology vanishings and work in this slightly broader setting; see \cite{AC} and \cite{Valerio+}.
In this way, in principle we loose the connection to determinantal representations of Chow forms, but note that we can always get back the very ample case by taking the embedding given by some sufficiently high multiple of $H$.
In particular, aCM bundles on $(X,H)$ stay aCM even on $(X,lH)$ for any $l\in\N$.
For Ulrich ones, this is not the case, but from Ulrich bundles on $(X,H)$ we can still construct Ulrich bundles on $(X,lH)$, which in general will have higher rank, see \cite[Prop. 5.4]{EisSch}.
From this point of view, it seems more natural to map $X$ to $\p^N$ with $|H|$ and then, after performing general projections, we get a finite and surjective ramified covering $f:X\to\p^n$, as granted by Noether's normalisation theorem, with the property that $\os_X(H)\cong f^*\osn(1)$.
Therefore, we are led to study varieties depending on the way they are presented as coverings of $\p^n$, i.e. depending on the structure of $f_*\os_X$.
We consider finite maps $f:X\to\p^n$ such that $X$ embeds as a divisor in $|\os_{P_m}(d)|$ for some $d>0$ where $P_m:=\p(\osn\oplus\osn(m))$ and $f=\pi|_X$ where $\pi:P_m\to \p^n$; we call them \textit{divisorial}.
Note that, we will see in \autoref{divisorialcover} that if $f$ is a divisorial covering then 
\[f_*\os_X\cong \osn\oplus\osn(-m)\oplus\dots\oplus\osn((1-d)m).\]
Remarkably, all cyclic coverings, the ones which on affine covers are described by the equation $t^d-f(\mathbf{x})$ in $Spec(\gk[x_1,\dots,x_n,t])$, are of this form.
Our first result is a generalisation of \cite[Thm. A]{Beadeterminantal}, using the same  strategy as in \cite[Prop. 8.1]{HanselkaKummer}:

\Te\label{1}(\autoref{Ulrichdivisorial},\autoref{ulrichmatrix''}) 
Let $f:X\to\p^n$ be an integral divisorial covering of degree $d$.
The following data are equivalent:
\begin{itemize}
    \item an aCM sheaf $\e$ of rank $r$ on $(X,H)$ such that $f_*\e\cong \oplus_{l=1}^{rd}\osn(\alpha_l)$
    \item a square matrix $A$ of order $rd$ whose entries are elements $a_{k,l}\in H^0(\os_{P_m}(1)\otimes\pi^*\osn(\alpha_k-\alpha_l))$ such that $p^r=det(A)$    
\end{itemize}
In particular, any such variety admits aCM and Ulrich sheaves (which correspond to the case $\alpha_l=0$ for all $l$).
\Ma

The last claim follows from the first one and the existence of matrix factorisations for arbitrary rings proved in \cite{BHU}.
Note that, existence of Ulrich bundles on double coverings of $\p^2$ have already been shown in \cite{SebastianTripathi}, for double coverings of $\p^n$ in \cite{KuNaPa} and on cyclic coverings of $\p^n$ in \cite{ParPin}.
However, in our result we get a better estimate on the minimal rank of such sheaves.
Moreover, notice that our methods work over any field.
As a further extension, we show some result of existence of Ulrich sheaves on some particular divisorial coverings of projectively normal subvarieties of $\p^N$ in \autoref{divgeneral} and, as a special example, we treat a class of Horikawa surfaces in \autoref{horik}.

More generally, given a finite surjective morphism $f:X\to\p^n$, by Horrock's theorem $\os_X$ has no intermediate cohomology if and only if $f_*\os_X$ is split.
Therefore, seen also the existence result of aCM sheaves by Faenzi--Pons-Llopis on aCM subvarieties of $\p^N$ \cite{FaenziPons_acm}, we can introduce the following modification of Eisenbud--Schreyer's conjecture, which seems more tractable at the moment: 

\medskip
\noindent\textbf{Question.} 
\textit{Does any finite covering $f:X\to\p^n$ with $f_*\os_X$ split admit an Ulrich sheaf with respect to $f^*\osn(1)$?
If not, what conditions can we ask on $f_*\os_X$ to grant existence of Ulrich sheaves?}
\medskip

The majority of this paper is devoted to the special situation of rank $2$ bundles on double coverings, with the aim of carefully analysing their geometry.
In this case, we can enhance our previous result to: 

\Te\label{2}(\autoref{rk2-section}, \autoref{section-rk2}, \autoref{23})
Suppose $char(\gk)\neq 2$ and $n\ge 3$. 
Let $f:X\to\p^n$ be a regular double covering with branch locus $B\subset \p^n$ of degree $2m$ and equation $b=0$.
If $Pic(X)\cong \Z H$ then for any $1\le\alpha,\beta\le m$ the following data are equivalent:
\begin{enumerate}
\item non-split rank $2$ aCM sheaf $\e$ on $X$ with $f_*\e\cong \osn\oplus\osn(\alpha-m)\oplus\osn(\beta-m)\oplus\osn(\alpha+\beta-2m)$
\item $Y\subset X$ mapped by $f$ isomorphically onto a complete intersection of two hypersurfaces of degrees $\alpha,\beta$
\item $b=p_\alpha q_{\alpha}+p_{\beta} q_{\beta}+p_m^2$ with $p_l\in H^0(\p^n,\osn(l))$ and $\{p_\alpha=0=p_\beta\}$ defines a complete intersection
\end{enumerate}
Moreover, Ulrich sheaves correspond exactly to the case $\alpha=\beta=m$.
\Ma

Note that from the last item we can easily construct a matrix as required in the previous theorem; see \autoref{esempio}.
Moreover, the complete intersection $p_\alpha=0=p_\beta$ will be the image of some $Y$ as in the second item.
The above $Y$ is actually the zero locus of a section of a unique $\e$, fact which translates into the exact sequence à la Hartshorne--Serre
\ses{\os_X}{\e}{\id_Y(\gamma)}
This link will be our main tool to analyse the existence of the sheaves $\e$ and some of their properties.
For example, we prove that the cones of Nef and effective divisors of $\p(\e)$ coincide and we study their second extremal ray \autoref{+} as well as show that $\e$ is fixed by the involution of $X$ if and only if $Y$ can be chosen inside the ramification divisor of $f$ \autoref{fisso'}.

In case $m=1$ then $X$ is a quadric; hence, we will skip it, even though our proofs would work equally well.
A dimensional count on the parameter space of the possible polynomials $b$ shows that for $n=3$ and $m>4$ or $n\ge 4$ and $m>1$ the general such $X$ cannot admit the presence of a rank $2$ Ulrich sheaf.
Since existence of Ulrich bundles of rank $2$ on complex double covering of $\p^2$ has already been proved in \cite{SebastianTripathi}, see also \cite[§5]{ParPin} for other cyclic covers of $\p^2$, we will focus on $3$-folds.
The cases in which we have hope for existence of rank $2$ Ulrich bundles on a general such $X$ are $m=2,3,4$, i.e. when $X$ is Fano or Calabi-Yau.
Over $\C$, the works \cite{ArrondoCosta}, \cite{madonna2001}, \cite{BF}, \cite{CH}, \cite{Kuznetsov_instanton}, \cite{Faenzi11}, \cite{Casnati_Faenzi_Malaspina_dP6_3}, \cite{CasnatiFilipMalaspina},\cite{Casnati_Faenzi_Malaspina_dP6_2} and \cite{Bea} ,  studied aCM and Ulrich bundles on Fano $3$-folds, with respect to the primitive submultiple of the anti-canonical bundle, facing the rank $2$ case.
As a result, besides existence of such sheaves, we have a classification of the possible Chern classes and information on the Hilbert schemes of zero loci of sections of those sheaves.
Moreover, in \cite{LMS}, \cite{Cho_Kim_Lee_ic22} and \cite{LeePark_UlrichV5} moduli spaces of Ulrich bundles of arbitrary rank $\ge2$ are studied using Bridgeland stability conditions on the Kuznetsov component of del Pezzo $3$-folds of degrees, respectively, $3,4,5$.
In the works \cite{CFK3} and \cite{CFK3.1}, the authors give a uniform argument for the classification of Ulrich bundles on Fano $3$-folds of index $i_X\ge 2$ or $i_X=1$ and cyclic Picard group, under the assumption that the primitive multiple of the anticanonical bundle is very ample.
Due to this last condition, some Fano $3$-folds, see \cite{IskPro} or \cite{fanography} for a classification, are not included in the previous works, e.g. double coverings of $\p^3$.
Nevertheless, some rank $2$ aCM bundles could be constructed from \autoref{2} and the present literature: for example, for $\alpha=1$ and $\beta=1,2$ those curves $Y$ are lines or conics whose existence on Fano $3$-folds is already known, see \cite[Chap. III]{Isk79} or \cite[§2]{KuzProShr}, while, when $\alpha=\beta=m=2$ those are elliptic curves of degree $4$ studied in \cite[§4]{Voisinqds}.
Both Gieseker and Bridgeland (on the Kuznetsov component) moduli space containing rank $2$ Ulrich bundles on the quartic double solid has been described in \cite[Appendix]{FLZ}.

In our work we give a uniform treatment for all $1\le\alpha,\beta\le m$ showing, for a general double covering $X$, not only the existence of such curves mapping isomorphically to complete intersections on $\p^3$ but also compute the cohomology of their normal bundles.
From that we can also study the regularity and dimension of the moduli spaces of the corresponding aCM sheaves.

\Te(\autoref{rk2esistenza}, \autoref{dimensioniacm})\label{3}
Suppose $char(\gk)\ne 2$.
Let $f:X\to\p^3$ be a general double covering with branch locus $B\subset \p^3$ of degree $2m=4,6,8$.
Then, for all $1\le\alpha,\beta\le m$ we can find the corresponding rank $2$ aCM and Ulrich sheaves $\e$.
Moreover, when they are stable (and $\gk=\overline{\gk}$ in case $m=4$), they belong to a generically regular component of their moduli spaces which has the expected dimension. 
\Ma
We can actually remove the assumption "$X$ general" in case $\alpha\ne \beta$ or $\alpha=1=\beta$, see \autoref{rk2esistenza}.
Remarkably, for $m=4$, i.e. the Calabi--Yau case, anytime the bundles we found are simple they are also \textit{spherical}, see \autoref{sferici}.
Such objects are central in the study of the derived category of $X$, see \cite[§8]{Huybrechts}.

To conclude, putting together the literature and the results of this work, in case $\gk=\C$ the only smooth Fano $3$-folds of index $2$ where no Ulrich bundles are known, at least to us, is the case \cite[\href{https://www.fanography.info/1-11}{1-11}]{fanography}, where the primitive submultiple of the anticanonical bundle has a base point.
On the other hand, among the ones with index $1$ and cyclic Picard group the minimum rank of an Ulrich bundle is always $2$\footnote{To be precise, \autoref{3} only takes into account general sextic double solids, but we believe it is an artefact of our proof} except (possibly) for \cite[\href{https://www.fanography.info/1-2}{1-2 b)}, \href{https://www.fanography.info/1-5}{1-5 b)}]{fanography}, where by \autoref{i1} the minimum can be shown to be at most $4$.
\medskip

\noindent\textbf{Structure of the paper.}
In Section $2$ we review some material on aCM and Ulrich sheaves while in Section $3$ we introduce the finite morphism, \textit{divisorial coverings}, we will be interested in.
The proof of \autoref{1} is in Section $4$, where also a new bound on the minimal rank of Ulrich sheaves on cyclic covers is given.
Section $5.1$ contains the proof of \autoref{2}. 
The rest of section $5$ is devoted to studying the action of the covering involution on rank $2$ aCM bundles and some properties of Ulrich ones, like syzygy bundles and cones of divisors of their projectification.
In Section $6$ we focus on $3$-folds and prove \autoref{3} relying on \textit{Macaulay2} codes in the appendix.
Finally, in Section $7$ we turn on some general application of our results to the existence of Ulrich bundles on varieties which can be constructed by pullback along divisorial coverings of $\p^n$.

\medskip
\noindent\textbf{Acknowledgements.}
This work is a generalisation of part of my PhD thesis at Tor Vergata University of Rome, see \cite{tesi}.
I thank my supervisor Ciro Ciliberto for his guidance and all the observation which improved the mathematical content as well as the readability of both the thesis and this work.
Moreover, I thank Nelson Alvarado, Fabio Bernasconi, Valerio Buttinelli, Davide Gori, Angelo Lopez and Antonio Rapagnetta for many useful conversations.
Part of this work have been done during my 3 months stay in University Paul Sabatier in Toulouse, so I thank Thomas Dedieu, Fulvio Gesmundo and Laurent Manivel for their hospitality and many useful conversations.
This work was supported by the MIUR Excellence Department Project MatMod@TOV awarded to the Department of Mathematics of the University of Rome Tor Vergata.
The author was also partially supported by GNSAGA of the Istituto Nazionale di Alta Matematica “F. Severi”.

\medskip
\noindent\textbf{Notation.}
Fix a field $\gk$, which for simplicity will be infinite.
If not otherwise stated, $(X,H)$ will always be a polarised $n$-dimensional variety over $\gk$, that is a proper, integral scheme of dimension $n$ over $\gk$ with an ample and base-point free Cartier divisor $H$.
We denote $\os_X(iH)$ simply by $\os_X(i)$ and tensorisation by $\os_X(i)$ simply by $(i)$.
As usual, $ext^i(-,-):=dim_\gk Ext^i(-,-)$ and $h^i(X,-):= ext^i(\os_X,-)$.
We denote by $\omega_X$ the dualising sheaf of $X$.
The symbol $\sim$ stands for linear equivalence.

%% file: generali.tex
\section{Ulrich and aCM sheaves}
In this section we recall definitions and basic properties of aCM and Ulrich sheaves.
\De\label{defUlrich}
A coherent sheaf $\e$ on a polarised variety $(X,H)$ of dimension $n$ is said: 
\begin{itemize}
    \item \textbf{initialised} if 
    \[h^0(X,\e(-1))=0 \quad \text{but} \quad h^0(X,\e)\neq 0;\]
    \item \textbf{has no intermediate cohomology} if
    \[h^j(X,\e(i))=0 \quad \text{for all} \quad i\in\Z,\; 1\leq j\leq n-1;\]
    \item \textbf{aCM} if it has no intermediate cohomology and for all $x\in X$ we have $depth(\e_x)=dim(\os_{X,x})$;
    \item \textbf{Ulrich} if 
    \[h^j(X,\e(i))=0 \quad \text{for}\quad 0\leq j\leq n,\, -n\leq i\leq -1.\]
\end{itemize}
\Ne

We can rephrase the aCM condition in purely cohomological terms.
We start with a lemma.

\Le\label{CM}
Let $\f$ be a coherent sheaf on a polarised variety $(X,H)$ such that $h^i(X,\f(-k))=0$ for all $i<n=dim(X)$ and $k$ big enough, then $depth_X(\f_x)=dim(\os_{X,x})$.
If, moreover, $X$ is regular, then $\f$ is locally free.
\ma
\begin{proof}
Let us embed $X$ in $\p^N$ using some multiple of $H$ and set $c:=N-n$ as the co-dimension of $X$.
We claim $Supp(\f)=X$.
Otherwise, if the support would have dimension $s<n$ then by our assumption and Serre duality, see \cite[Chap. 1 (1.3)]{AltKle_duality}, for $k\gg0$ we would have
\[0=h^s(X,\f(-k))=hom(\f(-k),\omega_X)=h^0(\Hom(\f(-k),\omega_X))=h^0(\Hom(\f,\omega_X)(k)).\]
But this contradicts the fact that $\Hom(\f,\omega_X)(k)$ is globally generated, being $H$ ample and $k\gg0$.
Therefore, for all $x\in X$ we have $\f_x\neq 0$ hence $dim(\os_{X,x})\ge depth_X(\f_x)$.

For the converse inequality it is enough to verify that $\ext^i(\f,\omega_{\p^N})=0$ for $i>c$: in such a case by \cite[Prop. 1.1.6 ii)]{HuyLeh} $\f$ would satisfy condition $S_{N-c,c}$ hence for any $x\in X$
\[depth_{\p^N}(\f_x)\ge min\{N-c, dim(\os_{\p^N,x})-c\}=dim(\os_{\p^N,x})-c=dim(\os_{X,x})\]
and to conclude we note that $depth_{\p^N}(\f_x)=depth_{X}(\f_x)$ by \cite[Rmk. 5.7.3 vi)]{EGA4.2}.
For $k\gg0$ the sheaf $\ext^i(\f,\omega_{\p^N})(k)$ is globally generated, being $H$ ample.
However, by \cite[Prop. III 6.9]{Har}, for $k\gg0$ we have 
    \[h^0(\p^N,\ext^i(\f,\omega_{\p^N})(k))=h^0(\p^N,\ext^i(\f(-k),\omega_{\p^N}))= ext^{i}(\f(-k),\omega_{\p^N})= h^{N-i}(\f(-k))\]
by Serre duality.
Finally, for $i>c$ we get $N-i<N-c=n$, so this last dimension is $0$ by assumption and we must have $\ext^i(\f,\omega_{\p^N})=0$ as desired.
    
     If $X$ is regular, then we can apply the Auslander-Buchsbaum formula, see \cite[Thm. 19.9]{EisComalg1}, and deduce $pd(\f_x)=0$, hence $\f$ is locally free. 
\end{proof}

As a corollary we deduce the well-known fact that aCM sheaves, and in particular Ulrich ones, are locally free when $X$ is regular.
The following is an extension of \cite[Prop. 1.6]{AnconaOttaviani} to the non-regular case.

\Prop\label{0tor_ext}
If $(X,H)$ is a polarised variety, any sheaf without intermediate cohomology $\e$ is an extension of a $0$-dimensional sheaf and an aCM sheaf.
If $X$ is regular then the extension is a direct sum.
\One
\begin{proof}
We have the exact sequence 
\sesl{T}{\e}{\f}{t}
where $T$ is the $0$-dimensional torsion subsheaf, see \cite[Def. 1.1.4]{HuyLeh}, in particular $\f$ will have no $0$-dimensional torsion.
   From \eqref{t} we see that $\f$ has no intermediate cohomology, since $T$ has no cohomology in degree higher than $0$.
   To conclude, we want to apply $\autoref{CM}$ to $\f$, so we need to verify that $h^0(\f(-i))=0$ for $i$ big enough.
Indeed, consider $s\in H^0(\f)$ and set $S:=Supp(s)$, then we have $dim(S)\ge 1$ since $\f$ has no $0$-dimensional torsion.
Being $|H|$ globally generated, we can choose some effective $D\in |H|$ such that $s|_D\neq 0$ and, by \cite[Lem. 1.1.12]{HuyLeh}, such that the sequence
\ses{\f(-H)}{\f}{\f|_D}
is exact.
Since the morphism $H^0(\f)\to H^0(\f|_D)$ is not the zero map, being $s|_D$ in its image, we get $h^0(\f(-H))<h^0(\f)$.
If $h^0(\f(-H))\ne 0$ we repeat the same operation, until $h^0(\f(-lH))=0$.
At this point, we can apply \autoref{CM}, which implies $\f$ is aCM.
If $X$ is regular then $\f$ is locally free by \autoref{CM} so that, by dimensional reasons, we have $ext^1(\f,T)=ext^1(\os_X,\f\du\otimes T)=h^1(\f\du\otimes T)=0.$
\end{proof}

\Co
A sheaf $\e$ without intermediate cohomology is aCM if and only if $h^0(\e(-l))=0$ for $l\gg0$. 
\io
\begin{proof}
If $\e$ is aCM then $depth(\e_x)=dim(\os_{X,x})$ hence $\e$ cannot have torsion subsheaves.
But then \autoref{0tor_ext} gives $T=0$ and $\e\cong\f$ is aCM.
    For the converse, taking cohomology in \eqref{t}, due to the fact that $dim(T)=0$ we get $h^0(T)=h^0(T(-l))\le h^0(\e(-l))=0$, hence $T$ must be the zero sheaf.
\end{proof}
Clearly, any extension of aCM (Ulrich) sheaves is aCM (Ulrich).
The notion of being aCM (Ulrich) has a nice functorial behaviour under pushforward for finite morphisms of polarised pairs, since those maps preserve cohomology.
Having an ample and globally generated line bundle $\os_X(H)$ is equivalent to the existence of some finite surjective map $f:X\to\p^n$ such that $f^*\osn(1)\cong \os_X(H)$, see Noether's normalisation theorem \cite[Thm. 13.89]{GW1}.
In practice, the complete linear system $|\os_X(1)|$ gives a finite morphism to some $\p^N$ and by composing with projections away from general points we get the above $f$.
Therefore, we will simply call $f$-aCM (Ulrich) an aCM (Ulrich) sheaf on $(X,f^*\osn(1))$.

\Le\label{proiezione}
Let $(X,H), (X',H')$ be two polarised varieties of dimension $n$.
Suppose that $g:X'\to X$ is a finite morphism such that $g^*H\sim H'$ and that $\e$ is a coherent sheaf on $X'$.
Then, $\e$ is $H'$-aCM (Ulrich) if and only if $g_*\e$ is $H$-aCM (Ulrich).
\ma
\begin{proof}
    Since the map $g$ is finite it preserves cohomology, see \cite[Cor. 22.6 1)]{GW2} hence, using \textit{projection formula} \cite[Thm. 22.81]{GW2}, we have 
    \[h^j(X',\e(iH'))=h^j(X',\e\otimes g^*\os_X(iH))=h^j(X,g_*(\e\otimes g^*\os_X(iH)))=h^j(X,g_*(\e)\otimes\os_X(iH)).\]
\end{proof}

It is a classical result of Horrocks, see \cite{Horrocks}, that aCM sheaves on $(\p^n,\osn(1))$ are all direct sums of line bundles, moreover, it is easily shown that all Ulrich bundles on such varieties are direct sums of $\osn$.
Therefore, applying \autoref{proiezione} we conclude the following.

\Co\label{noether}
A coherent sheaf $\e$ on $(X,\os_X(1))$ is aCM (Ulrich) if and only if $f_*\e$ is split (trivial) for some, equivalently any, finite morphism $f:X\to \p^n$ such that $f^*\osn(1)\cong \os_X(1)$.
\io

We immediately recover the known fact that Ulrich sheaves are initialised aCM sheaves with the maximum possible number of global sections, see \cite[Prop. 2.1]{EisSch}.

Ulrich sheaves have many remarkable properties.
First, we state some equivalent definitions that can be found in the literature.
A proof when $\os_X(1)$ is very ample can be found in \cite[Thm. 2.1]{EisSch}, \cite[Thm. 2.3]{Bea}, while for the general case we refer to \cite[Thm. 1.4]{AC}, see also \cite[§1.1]{tesi}.

\Te\label{definizioneUlrich}
Let $(X,H)$ be a polarised, $n$-dimensional variety over $\gk=\overline{\gk}$ and $\e$ a coherent sheaf on $X$.
The following are equivalent:
\begin{enumerate}[i)]
    \item $\e$ is an Ulrich sheaf for $(X,\os_X(H))$
    \item $h^j(X,\e(-j))=0$ for $1\leq j\leq n$ and $h^{j}(X,\e(-j-1))=0$ for $0\leq j\leq n-1$
    \item Denote $\phi:X\to\p^N$ the finite morphism given by the complete linear system $|H|$ and define $c:=N-n$, then we have a \textbf{linear resolution} 
    \[0\to \os_{\p^N}(-c)^{r_c}\to \dots \to \os_{\p^N}(-1)^{r_1}\to \os_{\p^N}^{r_0}\to \phi_*\e\to 0\]
    \item $\e$ is aCM, initialised and $h^n(X,\e(-n))=0$, equivalently its \textbf{cohomology table} reads: 
    \[ h^j(X,\e(i))=0\quad \text{if} \quad \begin{cases}
        j=0, \; i<0 \\
        1\leq j\leq n-1, \; i\in \Z \\
        j=n, \; i\geq -n.\\
        \end{cases}\]
\end{enumerate}
\Ma

Here we list some interesting properties of Ulrich sheaves.
Again, the proofs are similar to the very ample case in \cite{CH} and can be found in \cite[§1]{tesi}.

\Te\label{corollario}
If $\e$ is a rank $r$ Ulrich sheaf on an $n$-dimensional polarised variety $(X,H)$ with $d:=H^n$, then 
\begin{itemize}
    \item $\e$ is globally generated
    \item the Hilbert polynomial of $\e$ is $P(\e)=rd\binom{t+n}{n}$ and $h^0(X,\e)=rd$
    \item $\e$ is Gieseker-semistable, all its Jordan-Holder factors are Ulrich sheaves and if $\e$ is Gieseker-stable then it is also slope-stable 
    \item if $X$ is regular then $c_1(\e)\cdot H^{n-1}=\frac{r((n+1)H+K_X)H^{n-1}}{2}.$
\end{itemize}
\Ma

\section{Divisorial and cyclic coverings}\label{s_divisorial}
In this work, we prefer to consider our projective variety $X$ as a surjective finite ramified covering of $\p^n$ and we are interested in a particular class, which more resemble divisors in $\p^n$ and contains al the \textit{cyclic} ones.
In this section, we introduce these finite morphism and fix the notation used from this point onwards.
Set $P_m:=\p(\osn\oplus\osn(m))$, where $\p(-):=Proj(Sym(-))$ stands for the projectivisation parametrising $1$-dimensional quotients.
We denote $\pi:P_m\to\p^n$ the standard projection and $\os_{P_m}(1)$ the relative hyperplane bundle.
Finally, $H_0:=\p(\osn)\subset P_m$ is a copy of $\p^n$.

\De\label{defdiv}
Take a finite surjective map $f:X\to\p^n$ of degree $d$.
We will call $f$ (and by extension $X$) \textbf{divisorial} if there is some $m>0$ such that $f=\pi\circ i$ where $i:X\hookrightarrow P_m$ is the embedding of a divisor in $|\os_{P_m}(d)|$ not intersecting $H_0$.
We define $\os_X(1):=f^*\osn(1)$.
\Ne
Note that, by the definition of $H_0$, we have $\os_{P_m}(1)|_{H_0}\cong \os_{H_0}$, hence any divisor $\Delta\in|\os_{P_m}(1)|$ intersecting $H_0$ must contain it.
The same applies to $\Delta\in|\os_{P_m}(d)|$ for any $d>0$ hence, any integral $X\in |\os_{P_m}(d)|$ gives a divisorial covering, since the condition of not intersecting $H_0$ is automatic and implies $\pi|_\Delta$ is finite. 

\Oss
Removing $H_0$ from $\p_m$ we obtain the total space of the line bundle $\osn(m)$, see \cite[Chapter 8 Section 4]{EGAII}.
Therefore, asking $X$ to embed in this total space would result in a definition equivalent to the one given.
Yet in another way, the morphism given by $|\os_{\p_m}(1)|$ contracts $H_0$ and maps $\p_m$ to the cone over the $m$-th Veronese variety of $\p^{n}$, which is nothing else than the weighted projective space $\p(1^{n+1},m)$.
Hence, divisorial coverings are divisors in $\p(1^{n+1},m)$ that do not pass through the unique singular point of this space, the vertex.
\one

For future use, we state some properties of the varieties we are considering.
By \cite[Prop. 22.86]{GW2} and derived projection formula \cite[Porp. 22.84]{GW2}, for any $\g$ vector bundle on $\p^n$ we have 
\begin{equation}\label{projection}
    R^i\pi_*(\os_{P_m}(l)\otimes\pi^*\g)\cong R^i\pi_*(\os_{P_m}(l))\otimes\g=\begin{cases}
        Sym^l(\osn\oplus\osn(m))\otimes\g \hfill i=0, l\ge 0 \\
        \left(Sym^{-2-l}(\osn\oplus\osn(m))\right)\du\otimes\g(-m) \qquad \hfill i=1, l\le -2 \\
        0 \hfill \text{otherwise}\\
    \end{cases}
\end{equation}
Those help us to study divisorial coverings.
\Le\label{divisorialcover}
    Let $f:X\to \p^n$ be a divisorial covering of degree $d\ge 1$ inside $P_m=\p(\osn\oplus\osn(m))$ then 
    \begin{enumerate}[1)] 
        \item $X$ is Gorenstein, $\os_{P_m}(d)|_X\cong \os_X(dm)$ and $\omega_X\cong \os_X(m(d-1)-n-1))$ 
        \item $f$ is flat and $f_*\os_X\cong \bigoplus_{i=0}^{d-1}\osn(-im)$ 
        \item $h^0(\os_X)=1$, in particular $X$ is connected, and $\os_X$ is aCM on $(X,\os_X(1))$.
    \end{enumerate}
    \ma
    \begin{proof}
        \gel{1}
        Since $X$ is a Cartier divisor in $P_m$, which is regular, it must be Gorenstein.
        By definition, we have $H_0\cap X=\emptyset$ hence $\os_{P_m}(H_0)|_X\cong \os_X$.
        We know that $\os_{P_m}(H_{0})\cong \os_{P_m}(1)\otimes\pi^*\osn(-m)$, see for example \cite[Prop. 4.6.2]{Sernesi}, therefore $\os_{P_m}(1)|_X\cong f^*\osn(m)$ and the first claim follows for any $d$ taking tensor products.  
        
        From \cite[Thm. 22.86 4)]{GW2} we know $\omega_{P_m/\p^n}  \cong\os_{P_m}(-2)\otimes\pi^*(\osn(m))$.
        Since $P$ is smooth then the sequence of relative differentials for the map $\pi$ is exact and taking determinants gives \[\omega_{{P_m}}  \cong\omega_{P_m/\p^n}\otimes\pi^*\omega_{\p^n}\cong\os_{P_m}(-2)\otimes\pi^*(\osn(m-n-1)).\]
         By the adjunction formula, \cite[Cor. 25.130 2)]{GW2}, we have \[\omega_X\cong \omega_{P_m}(d)|_X\cong \os_{P_m}(d-2)\otimes\pi^*(\osn(m-n-1))|_X\cong \os_X(m(d-1)-n-1)\] where we also used $1)$.

\gel{2}
By definition, the sequence defining $X$ in $P$ is
\ses{\os_{P_m}(-d)}{\os_{P_m}}{\os_X}
If we apply $\pi_*$ to it, using \eqref{projection} we obtain
\sesl{\osn}{f_*\os_X}{\bigoplus_{i=1}^{d-1}\osn(-im)}{d}
Note that $f$ is affine and $f_*\os_X$ is locally free so $f$ must also be flat, see \cite[Prop. 12.19]{GW1}.
Moreover, \eqref{d} splits being
$Ext^1\left(\bigoplus_{i=1}^{d-1}\osn(-im),\osn\right)\cong H^1\left(\bigoplus_{i=1}^{d-1}\osn(im)\right)=0.$

\gel{3}
The map $f$ is finite hence preserves cohomology then, by $2)$ we get $h^0(X,\os_X)=h^0(\p^n,\osn)=1$ hence $X$ is connected.
$\os_X$ is aCM by $2)$ and \autoref{noether}.
\end{proof}

Conversely, we can identify a divisorial covering just by looking at $f_*\os_X$.
\Le\label{1div}
A finite morphism of degree $d$ is divisorial if and only if $f_*\os_X\cong \bigoplus_{i=0}^{d-1}\osn(-im)$ for some $m>0$.
\ma
\begin{proof}
For the first claim, by the universal property of the symmetric algebra, the morphism $\osn(-m)\to f_*\os_X$ gives us a surjection $Sym(\osn(-m))\twoheadrightarrow f_*\os_X$ of $\osn$-algebras, which results in a closed embedding 
\[X\subset \mathrm{S}pec(Sym(\osn(-m)))\subset \p(\osn\oplus\osn(m))\]
avoiding $H_0=\p(\osn)$.
\end{proof}

We can parametrise divisorial coverings by polynomials as follows.
Call $t$ the section of $H^0(P_m,\os_{P_m}(1))\cong H^0(\p^n,\osn\oplus\osn(m))$ corresponding to $H^0(\p^n,\osn)$.
We have
\[H^0(P_m,\os_{P_m}(d))\cong H^0(\p^n,\pi_*\os_{P_m}(d))=H^0\left(\p^n,\bigoplus_{j=0}^dt^j\otimes\osn\left((d-j)m\right)\right)=\bigoplus_{j=0}^dt^jH^0(\p^n,\osn((d-j)m))\]
Therefore, if $x_1,\dots , x_n$ are coordinates on $\p^n$, we get that $H^0(P_m,\os_{P_m}(d))$ is identified with the homogeneous polynomials of degree $dm$ in the ring $\gk[x_0,\dots, x_n,t]$, where $x_i$ have degree $1$ and $t$ has degree $m$.
Hence, to any divisorial covering $X$ we can associate some polynomial, which we call an \textbf{equation} for $X$ (clearly well defined only up to a constant). 
Among such coverings, there is a special class.

\De\label{defcyclic}
A degree $d$ divisorial covering $f:X\to\p^n$ is \textbf{cyclic} if $X$ has equation in $P_m$ of the form $t^d-b$, with $b\in H^0(\p^n,\osn(dm))$.
\Ne

These equations give us a parameter space for divisorial coverings, through them we can define a notion of "generality".
\De\label{generaldivisorial}
Fixed $n,d,m$, we say that some property $\mathscr{P}$ \textbf{holds for the general divisorial covering} if: the set $\{p\in H^0(P_m,\os_{P_m}(d)) \;|\; \mathscr{P} \; \text{holds on the divisorial covering determined by}\; p\}$ contains an open subset for the Zariski topology.
\noindent Fixed $n,d,m$, we say that some property \textbf{holds for the general cyclic covering} if: the set $\{b\in H^0(\p^n,\osn(md)) \;|\; \mathscr{P} \; \text{holds on the cyclic covering determined by}\; b\}$ contains an open subset for the Zariski topology.
\Ne

Let us show some well-known features of cyclic covers.
The \textbf{ramification divisor} $R$ of a cyclic covering $f:X\to\p^n$ is the divisor given by $\{t=0\}$, where $t$ is seen as a section of $\os_{P_m}(1)$.
The divisor $B:=f(R)\cong R$ is called \textbf{branch locus} and $b$ is an equation for it in $\p^n$.
Note that $f^*B=dR$ as divisors.
A local computation shows that $R$ is exactly the locus on which the morphism $f$ is ramified, if $char(\gk)$ is coprime with the degree $d$ of $f$.

A finite morphism $f:X\to\p^n$ of degree $2$ will be called \textbf{double covering}.

\Le\label{2div}
Any flat double covering $f:X\to\p^n$ with $h^0(\os_X)=1$ is divisorial.
Moreover, if $char(\gk)\neq 2$ any divisorial double covering is cyclic.
\ma
\begin{proof}
By assumption we have $f_*\os_X$ locally free of rank $2$ and an exact sequence
     \ses{\osn}{f_*\os_X}{\f}
     Since the left morphism never drops rank, we have $\f\cong \osn(-m)$ for some $m\in\Z$ and hence the sequence splits.
     Being $1=h^0(\os_X)=h^0(f_*\os_X)$ then $m>0$, hence we can apply \autoref{1div}.

     For the second part, just complete the square in its equation, which is allowed since $2$ in invertible.
\end{proof}

On a cyclic double covering, there is an action of $\Z/2$, which can be explicitly constructed as follows.
On the weighted projective space $\p(1^{n+1},m)$ we can define an automorphism by sending $x_i\mapsto x_i$ and $t\mapsto -t$ and the form of the equation of $X\subset \p(1^{n+1},m)$ implies that $X$ is mapped to itself.

Finally, let us recall that often the Picard group of those divisorial coverings is cyclic and generated by $\os_X(1)$.
\Le\label{Piccover}
Suppose $f:X\to\p^n$ is a regular divisorial covering of degree $d$.
$Pic(X)$ is generated by $\os_X(1)$ if at least one of the following holds:
\begin{itemize}
    \item $f$ is cyclic, $char(\gk)$ is coprime with $d$ and $n\geq 4$
    \item $\gk=\C$ and $n>d$.
\end{itemize}
\ma
\begin{proof}
In the first case, it follows from \cite[XII Corollaire 3.7]{SGA2} that $Pic(B)\cong Pic(\p^n)$, where $B$ is the branch divisor of $f$.
From the sequence
\ses{\os_X(-R)}{\os_X}{\os_R}
together with the fact that $(X,\os_X(1))$ is aCM, see \autoref{divisorialcover}, and $n\ge 4$ we conclude that $H^j(R,\os_R(-i))=0$ for $j=1,2$ and $i>0$.
Together with the fact that $B$ is ample on $X$ and $B\cong R$, this implies that we can apply \cite[XII Corollaire 3.6]{SGA2} and get $Pic(R)\cong Pic(X)$ hence the claim.
In the second case, we are in a position to apply \cite[Prop. 3.1]{Laz}. 
\end{proof}

\section{aCM and Ulrich sheaves on divisorial coverings}

For divisors in $\p^n$ we have a dictionary between aCM and Ulrich sheaves and determinantal presentations of their equations.
We extend those result to divisorial coverings exploiting their similarity to divisors in $\p^n$.
As a result we get \autoref{1}.
Moreover, for cyclic coverings we can also improve the current estimate for the minimal rank of Ulrich sheaves.

\subsection{Matrix representation}
Suppose $f:X\to\p^n$ is the embedding of a divisor.
We know that a sheaf $\e$ on $X$ is aCM if and only if we have a linear resolution 
\ses{\bigoplus_{i}^l\osn(\alpha_i)}{\bigoplus_{i}^l\osn(\beta_i)}{\e}
see \cite[Thm. A]{Beadeterminantal}.
Moreover, $\e$ is Ulrich if and only if we can choose $\alpha_i=-1$ and $\beta_i=0$ for all $i$, see \cite[Prop. 2.1]{EisSch}.
Here we want to derive a similar result in the case where $f$ realises $X$ as a divisorial covering of $\p^n$.
Since $X$ is a divisor in $P_m$, the aCM sheaf $\e$ has a locally free resolution of length $1$ on $P_m$, our goal is to give a precise description of it in terms of the splitting type of $\e$, that is, of $f_*\e$.

Historically, the link between aCM modules over hypersurface local rings and determinantal representations of their equations has already been pointed out in \cite[§6]{Eis-resolutions}.
The next result for Ulrich sheaves on cyclic coverings of $\p^n$ is essentially contained in \cite[Prop. 8.1, Remark 8.2]{HanselkaKummer}.
Here we generalise this result to aCM sheaves on divisorial coverings.

\Te\label{Ulrichdivisorial}
Let $f:X\to\p^n$ be an integral divisorial covering of degree $d$ such that $i:X\hookrightarrow \p_m$ has equation $p=0$.
If there are $\alpha_l\in\Z$ for $l=1,\dots ,rd$ and a matrix $A$ of order $rd$ whose $(k,l)$-entry is $a_{k,l}\in H^0(\os_{P_m}(1)\otimes\pi^*\osn(\alpha_k-\alpha_l))$ such that $p^r=det(A)$, then on $X$ there is an $f$-aCM sheaf $\e$ of rank $r$ such that $f_*\e\cong \bigoplus_{l=1}^{rd}\osn(\alpha_l)$ and fits in
\begin{equation}\label{resol}
    0\to \left(\bigoplus_{l=1}^{rd}\pi^*\osn(\alpha_l)\right)\otimes\os_{P_m}(-1)\overset{A\cdot }{\longrightarrow} \bigoplus_{l=1}^{rd}\pi^*\osn(\alpha_l)\to i_*\e\to 0.
\end{equation}
Conversely, if $\e$ is an $f$-aCM sheaf of rank $r$ on $X$ such that $f_*\e\cong \bigoplus_{l=1}^{rd}\osn(\alpha_l)$, then we can find a matrix $A$ as above and a resolution as in \eqref{resol}.
\Ma
\begin{proof}
    For the first claim, a matrix $A$ as described gives a morphism of sheaves as in \eqref{resol}, which is injective and has cokernel supported on $X$ since $det(A)=p^r$ is not identically zero.
    By \autoref{noether}, to verify that $\e$ is $f$-aCM it is enough to show that $f_*\e= \pi_*i_*\e$ is split, which follows by applying $\pi_*$ to \eqref{resol} and recalling \eqref{projection}.
    Being $X$ integral, we have $c_1(\e)\sim rk(\e)[X]\sim d\cdot c_1(\os_{P_m}(1))rk(\e)$ but then, we get $rk(\e)=r$ from
    \[drk(\e)\cdot c_1(\os_{P_m}(1))\sim  c_1(\e)\sim c_1(\os_{P_m}^{rd})-c_1(\os_{P_m}(-1)^{rd})\sim rd\cdot c_1(\os_{P_m}(1)).\]
For the converse, suppose $f_*\e\cong \oplus_{l=1}^{rd}\osn(\alpha_l)$.
By adjunction between $\pi_*$ ad $\pi^*$, to the identity morphism of $f_*\e\cong \pi_*i_*\e$ there corresponds a morphism $\pi^*\oplus_{l=1}^{rd}\osn(\alpha_l)\cong \pi^*\pi_*i_*\e\xrightarrow{\phi} i_*\e$.
First, we will show that it is surjective.
To this end, it is enough to show that $\phi$ is surjective when restricted to each fiber of $\pi$.
Let $y\in\p^n$ be a closed point, $\pi_y:=\pi\inv(y)$ the fiber above $y$ and $j:\pi_y\hookrightarrow \p_m$ its inclusion.
Note that, being $f$ finite of degree $d$ and $\e$ of rank $r$, for any $y\in \p^n$ the sheaf $j^*i_*\e$ is supported on a finite length scheme and is actually a module of length $rd$ over $\gk$, in particular the dimension of its global sections is constantly equal to $rd$.
Moreover, $\e$ is flat over $\p^n$: this is a local property and $f$ is affine, then the claim follows from the fact that $f_*\e$ is locally free.
Therefore, by Grauert's theorem \cite[Thm. 23.140]{GW2}, $\pi_*i_*\e\otimes k(y)\cong H^0(j^*i_*\e)$.
By commutativity of 
\[\xymatrix{\pi_y \ar[r]^j \ar[d] & P_m \ar[d]^\pi\\
y \ar[r] & \p^n\\}\]
we have $j^*\pi^*\pi_*i_*\e\cong (\pi_*i_*\e\otimes k(y))\otimes\os_{\pi_y}$ hence, the morphism $\phi$ once restricted to $\pi_y$ is the evaluation map of the global sections of $j^*i_*\e$, which is surjective being this last sheaf supported on $0$-dimensional schemes.

We can now form the short exact sequence
\begin{equation}\label{resol'}
0\to\ck\to\bigoplus_{l=1}^{rd}\pi^*\osn(\alpha_l)\cong \pi^*f_*\e\xrightarrow{\phi}i_*\e\to 0
    \end{equation}
and we are left to show that $\ck\otimes\os_{P_m}(1)\cong \pi^*f_*\e\cong \pi^*\oplus_{l=1}^{rd}\osn(\alpha_l)$, so that the left morphism in \eqref{resol'} would be given by a matrix of the desired form.
Apply $\pi_*$ to \eqref{resol'}.
Using \eqref{projection} and $\pi_*i_*\e=f_*\e$ we get 
\begin{equation}\label{push'}
    0\to \pi_*\ck\to \pi_*\left(\bigoplus_{l=1}^{rd}\pi^*\osn(\alpha_l)\right)\xrightarrow{\pi_*\phi} f_*\e\to R^1\pi_*\ck\to 0
\end{equation}
We will show that $\pi_*(\phi)$ is an isomorphism.
Again by adjunction of $\pi_*$ and $\pi_*$, we have a composition 
\[\pi_*(i_*\e)\to \pi_*\pi^*\pi_*(i_*\e)\xrightarrow{\pi_*\phi} \pi_*(i_*\e)\]
which is the identity, in particular the second map is surjective.
Being $\pi_*\os_{\p_m}\cong \osn$, by projection formula we have $\pi_*\pi^*\pi_*i_*\e\cong \pi_*i_*\e$ hence the above maps have to be all isomorphisms, as desired.
We deduce by \eqref{push'} that $R^i\pi_*\ck=0$ for all $i= 0,1$.
Combining \cite[Cor. 23.114]{GW2} and \cite[Thm. 23.140 1)]{GW2}, the restrictions of $\ck$ to the fibers of $\pi$ have no cohomology hence have to be $\on{1}(-1)^{rd}$, being such fibers $\p^1$-s.
Therefore, the restrictions of $\ck\otimes\os_{P_m}(1)$ to these fibers are trivial.
But then, \cite[Prop. 25.1.11]{Vakil}\footnote{this is stated only for line bundles but the same proof works for locally free sheaves of any rank.} implies that $\ck\otimes\os_{P_m}(1)\cong \pi^*\pi_*(\ck\otimes\os_{P_m}(1))$.
To compute $\ck$, let us apply $\pi_*$ to \eqref{resol'} twisted by $\os_{P_m}(-1)$.
Observe that $\os_{P_m}(-1)|_X\cong \os_X(-m)$ by \autoref{divisorialcover} $1)$.
Moreover, $R^i\pi_*(\pi^*f_*\e\otimes\os_{P_m}(-1))\cong R^i\pi_*(\os_{P_m}(-1))\otimes f_*\e=0$ for all $i$ by \eqref{projection}.
Therefore, we conclude 
\[f_*(\e)(-m)\cong R^1\pi_*(\ck\otimes\os_{P_m}(-1))\cong R^1\pi_*(\os_{P_m}(-2)\otimes(\ck\otimes\os_{P_m}(1)))\cong \]
\[\cong R^1\pi_*(\os_{P_m}(-2)\otimes\pi^*\pi_*(\ck\otimes\os_{P_m}(1)))\cong \pi_*(\ck\otimes\os_{P_m}(1))(-m)\]
where in the last step we used again \eqref{projection}.
We conclude $\pi_*\ck\otimes\os_{P_m}(1)\cong f_*\e$ whence $\ck\cong \pi^*(f_*\e)\otimes\os_{P_m}(-1)$.
\end{proof}

As a corollary, we can prove that the locus parametrising divisorial coverings is constructible.
\Co\label{constructible}
Fixed $n,m,d,r$, the locus in $H^0(\os_{P_m}(d))$ corresponding to divisorial coverings having a rank $r$ aCM sheaf is constructible; if $r=1$ it is also irreducible.
\io
\begin{proof}
Consider $U:=\bigoplus_{l,k}H^0(\os_{P_m}(1)\otimes\pi^*\osn(\alpha_k-\alpha_l))$, the space parametrising square matrices of order $rd$ whose entry $a_{k,l}$ is in $H^0(\os_{P_m}(1)\otimes\pi^*\osn(\alpha_k-\alpha_l))$.
Taking the determinant gives an algebraic map $U\to H^0(\os_{P_m}(rd))$ whose image, call it $T$, is clearly irreducible, and is constructible by Chevalley's theorem \cite[Thm. 10.20]{GW1}.
If $r=1$ by \autoref{Ulrichdivisorial} we are done, since the intersection of $T$ and the open set of divisorial coverings is still constructible and irreducible.
Otherwise, we take the intersection of $T$ with the subvariety of $H^0(\os_{P_m}(rd))$ parametrising $r$-th powers of polynomials in $H^0(\os_{P_m}(d))$, which is again constructible.
\end{proof}

\subsection{Matrix factorizations and Ulrich sheaves on divisorial coverings}

We want to focus our attention on the matrix $A$ and on its relation to the polynomial $p$.
This will allow us to prove existence of aCM and Ulrich sheaves using the existence results on matrix factorisations.
We start by recalling an important definition.
In the following, $\gi_l$ is the identity matrix of order $l$.
\De
Let $\mathfrak{R}$ be a commutative unital ring and $\mathfrak{S}\subset \mathfrak{R}$ any subset.
A \textbf{matrix factorisation} of an element $p\in \mathfrak{R}$ with coefficients in $\mathfrak{S}$ is the datum of $A_1,\dots , A_l$ matrices with entries in $\mathfrak{S}$ such that $\prod_{i=1}^l A_i=p\gi_n$ for some $n\in \N$.
If all $A_i$ are square matrices of the same order then this order will be the \textbf{size} of the factorization.
\Ne
Note that if $A_i$ are square matrices then
\[\prod_i det(A_i)=det\left(\prod_iA_i\right)=det(p\gi_n)=p^n\]
then, when $p$ is indecomposable, any matrix in a matrix factorisation for $p$ produces a determinantal representation for some power of $p$.
From \cite[§6]{Eis-resolutions}, the existence of maximally Cohen-Macaulay modules on the quotient of a regular local ring by the ideal generated by a single element $p$ is equivalent to the existence of a matrix factorizations of $p$.
Moreover, the matrices appearing in the factorization of $x$ give a resolution of the above module.
We start by restating our previous results using the language of matrix factorizations; this is the analogous of \cite[Cor. 6.3]{Eis-resolutions}.

\Co\label{factorulrich}
Consider an integral divisorial covering $f:X\to\p^n$ of degree $d$ and with equation $p=0$ in $P_m$.
There is a matrix factorisation of size $rd$ of the polynomial $p$ such that at least one of the matrices has coefficients $a_{k,l}\in H^0(\pi^*\osn(\alpha_k-\alpha_l)\otimes\os_P(1))$ if and only if there is an aCM sheaf $\e$ on $X$.
\io
\begin{proof}
If there is a matrix factorization of $p$ as above, call $A$ the special matrix whose existence is granted by hypothesis.
Being $X$ integral, $p$ is indecomposable and $det(A)=p^r$, so we conclude by \autoref{Ulrichdivisorial}.
    
Conversely, if there is an aCM sheaf $\e$ then we have the resolution \eqref{resol}.
Since $\os_{\p_m}(-X)\cong \os_{\p_m}(-d)$, the matrix $p\gi_{rd}$ induces a morphism $\os_{\p_m}(-d)\otimes\pi^*f_*\e\to\pi^*f_*\e$, whose image goes to $0$ when composed with the map $\pi^*f_*\e\to i_*\e$. 
    Therefore, we have a commutative diagram as in
    \dia 
    0 \ar[r] & \os_{\p_m}(-d)\otimes\pi^*f_*\e \ar[r,"id"] \ar[d] & \os_{\p_m}(-d)\otimes\pi^*f_*\e \ar[d,"p\gi_{rd}"] & &  \\
    0 \ar[r] & \os_{\p_m}(-1)\otimes\pi^*f_*\e \ar[r,"A"] &  \pi^*f_*\e \ar[r] & i_*\e  \ar[r] & 0. \\
    \mma 
    We conclude that there is a matrix $A'$, determined by the just constructed left vertical map, of size $rd$ such that $AA'=p\gi_{rd}$, hence we have the desired matrix factorization.
\end{proof}

Now we prove that any divisorial covering admits not only aCM but even Ulrich sheaves.
Under stronger assumptions, the case of double covering has been proved in \cite{KuNaPa} and the cyclic case in \cite{ParPin}.
\Te\label{ulrichmatrix''}
Any integral divisorial covering $f:X\to\p^n$ admits an $f$-Ulrich sheaf, in particular this holds if $X$ is a cyclic cover.
More precisely, if $d=deg(f)$ and for some $s\geq 2$\footnote{when $s=1$ we must have $X$ reducible, actually a union of copies of $\p^n$} we have
\[p=\sum_{i=1}^s\prod_{j=1}^dp_{i,j},  \qquad p_{i,j}\in H^0(\os_{P_m}(1))\] 
then there are Ulrich sheaves of rank $d^{s-2}\cdot\varphi(d)$, where $\varphi$ is Euler's totient function.
If $\gk$ contains $d$-th roots of unity then the rank can be taken to be just $d^{s-2}$.
\Ma
\begin{proof}
Call $p\in H^0(\os_P(d))$ the equation of $X$.
From \autoref{factorulrich}, to get an Ulrich sheaf of rank $r$ it is enough to construct a matrix factorisation of $p$ of order $rd$ and with coefficients in $H^0(\os_P(1))$.
Note that $H^0(\osn(m))$ generates $\bigoplus_{i\in\N}H^0(\osn(im))$ as a $\gk$-algebra hence $H^0(\os_P(1))$ generates $\bigoplus_{i\in\N}H^0(\os_P(i))$ as a $\gk$-algebra by \eqref{projection}.
Therefore, we can always write $p$ in the above form for some $s$ and $p_{i,j}\in H^0(\os_P(1))$.
But then we just apply \cite[Lem. 1.6]{BHU}, or \cite[Lem. 1.5]{BHU} in the case where $\gk$ contains $d$-th roots of unity, and we get the existence and the claimed ranks.
\end{proof}

\Co\label{doppi}
If $\gk$ is algebraically closed then any flat double covering $f:X\to\p^n$ always has $f$-Ulrich sheaves.
\io
\begin{proof}
    Being $\gk$ algebraically closed we have $h^0(\os_X)=1$ hence such covering is divisorial by \autoref{2div}, then we apply the previous theorem.
\end{proof}

\Co\label{triple}
For $n\geq 4$ and $\gk=\C$, every finite morphism $f:X\to\p^n$ of degree $3$ with $X$ smooth has an $f$-Ulrich bundle.
\io
\begin{proof}
     This $X$ can always be embedded as a divisor in the total space of a line bundle over $\p^n$ by \cite[Prop. 3.2]{Laz} and hence is a divisorial covering by definition.
\end{proof}

We cannot hope for similar results for higher degree coverings because Lazarsfeld gave examples of degree $5$ covering of $\p^n$ for arbitrary high $n$ which are not divisorial, see \cite[Remark 3.5]{Laz}.
Still, the question of lower-dimensional triple coverings remains open.
\Do
Does every flat triple covering of $\p^2$ or $\p^3$ carry an Ulrich sheaf?
\da

\subsection{Finer rank estimate for cyclic coverings}

Next, we will specialise in the case of cyclic coverings.
If we are interested in the minimal rank of an Ulrich sheaf on a divisorial covering $f:X\to\p^n$ then, from the previous results, we need to search for matrix factorisation of minimal size for its equation $p$.
This minimal rank highly depends on the specific form of $p$, indeed, in the proof of \autoref{ulrichmatrix''} we exploited the fact that $p$ can be written as a sum of $s$ products of forms in $H^0(P,\os_P(1))$.
For cyclic coverings the equation of $X$ is of the form $p=t^d-b$ with $b\in H^0(\p^n,\osn(dm))$, so we can find a more suitable way to express $p$.
Our point is that, since one of the summands in the equation for $p$ in \autoref{ulrichmatrix''} is equal to $t^d$, we can further subtract a term $p_0^d$ without increasing the rank of the matrix needed.

\Prop\label{ulrichsommeprodotti}
Let $f:X\to\p^n$ be an integral cyclic covering of degree $d$ over a field  containing a $d$-th root of unity.
If there is some $s\ge2$ such that the equation $b$ of the branch locus of $f$ can be written as $b=p_0^d+\sum_{i=1}^{s-1}\prod_{j=1}^d p_{i,j}$, where $p_{i,j}\in H^0(\p^n,\osn(m))$, then there exists an $f$-Ulrich sheaf of rank $d^{s-2}$.
\One
\begin{proof}
Take a $d$-th root of unity $\zeta_d$.
Then $t^d-p_0^d=\prod_{j=1}^d(t-\zeta_d^jp_0)$ therefore $p=t^d-b$ can be written as a sum of $s$, instead of $s+1$, products of forms in $H^0(P_m,\os_{P_m}(1))$.
By \autoref{ulrichmatrix''} there is an Ulrich sheaf of rank $d^{s-2}$.
\end{proof}

As an example we consider the lowest possible values for $s$.
\Es\label{esempio}
For $s=2$ then we have $p=t^d-b=t^d-p_0^d-\prod_{j=1}^dp_j$ and there is a rank $1$ Ulrich sheaf on such on $X$.
For $d=2,3$, following the construction in \cite{BHU} we get that $p$ can be expressed, respectively, as the determinant of the following matrices
\begin{equation}
    \begin{pmatrix}\label{det21}
    t-p_0 & p_1\\
    p_2  &  t+p_0\\
\end{pmatrix} \qquad \quad \begin{pmatrix}
    t-p_0 & p_1 & 0 \\
    0 & t-\zeta_3p_0 & p_2 \\
    p_3  & 0 & t-\zeta_3^2p_0\\
\end{pmatrix}.
\end{equation}
In general, we take a matrix having $t-\zeta_d^ip_0$ on the diagonal and the polynomials $p_i$ on the above-diagonal and left down corner.

For $s=3$ and $d=2$ we have $b=p_0^d+\sum_{i=1}^{2}\prod_{j=1}^d p_{i,j}$ and there is a rank $2$ Ulrich sheaf on $X$.
Therefore, we can write $(t^d-b)^2$ as a determinant but even more, as in \cite[Cor. 2.4]{Bea}, this polynomial can be seen as a Pfaffian of some skew-symmetric matrix.
Let us illustrate this in case $d=2$. 
We have $t^2-b=t^2-p_0^2-p_1p_2-p_3p_4$ hence from the construction in \cite{BHU} we get $(t^2-b)^2$ equal to
\[det\begin{pmatrix}
    t-p_0 & 0 & p_1 & p_3 \\
    0 & t-p_0 & -p_4 & p_2 \\
    p_2  & -p_3 & t+p_0 & 0 \\
    p_4  & p_1 & 0 & t+p_0 \\
\end{pmatrix}=det\begin{pmatrix}
     0 & t-p_0 & -p_4 & p_2 \\
    -(t-p_0) & 0 & -p_1 & -p_3 \\
    p_4  & p_1 & 0 & t+p_0 \\
    -p_2  & p_3 & -(t+p_0) & 0 \\
\end{pmatrix}\]
where the second matrix is obtained simply by exchanging first with second and third with fourth row, then multiplying by $-1$ the second and fourth.
\Io

Let us recall the state of the art for double coverings.
\Oss
Suppose $d=2$.
The best upper bound for the minimal rank of an Ulrich bundle on a regular double covering of $\p^n$ was the one given in \cite[Lem. 6.2]{KuNaPa}, their result reads as follows.
If the equation of the branch locus is $b=\sum_{i=1}^{s-1} p_{i,1}p_{i,2}$ then there is an Ulrich sheaf of rank $2^{s-1}$ or $2^{s-2}$.
Our proof implies that, to have an Ulrich sheaf of rank $2^{s-2}$ is sufficient that $b=p_0^2+\sum_{i=1}^{s-1} p_{i,1}p_{i,2}$.
Although it is only a slight improvement, it will be optimal in the cases $n=3$ and $m=2,3,4$, see \autoref{polinomi}.
\one

%% file: double.tex
\section{Characterising rank $2$ aCM bundles}

Fix a double covering $f:X\to\p^n$ with $n\ge 3$, which is naturally polarised by $\os_X(1):=f^*\osn(1)$.
Our goal is to characterise rank $2$ aCM bundles, in particular Ulrich ones, on such varieties using Hartshorne--Serre correspondence, i.e. studying the zero loci of their sections, and describe their geometry.
This will prove \autoref{2}.
Moreover, for cyclic double coverings, we describe when such sheaves are fixed by the involution.
Finally, we investigate the positivity of the Ulrich ones.

The zero loci of sections of $\e$ will be denote $Y$, with the convention that $\os_Y(l)=\os_X(l)|_Y$.
However, if $f_*\os_X\cong \osn\oplus\osn(-m),$ then we can think $X\subset P_m$ and we have $\os_{P_m}(1)|_X\cong \os_X(m)$ by \autoref{divisorialcover}, so that $\os_{P_m}(1)|_Y\cong \os_Y(m)$.

\subsection{Rank $2$ aCM bundle through Hartshorn-Serre correspondence}\label{s5.1}

As usual, for $\alpha,\beta\in\Z_{>0}$ we call $Z\subset \p^n$ a \textbf{complete intersection of type} $\mathbf{(\alpha,\beta)}$ if we have a short exact sequence of the form
\sesl{\osn(-\alpha-\beta)}{\osn(-\alpha)\oplus\osn(-\beta)}{\id_Z}{ic}

First, we characterise aCM bundles for smooth cyclic double coverings with cyclic Picard group.
Note that $\os_X(i)$ are always aCM line bundles on any divisorial covering, hence also their direct sums are, but we are interested in indecomposable, possibly stable, ones.
Moreover, it is easily verified that $\os_X(i)$, are never Ulrich.
\Prop\label{rk2-section}
Over a field $\gk$, let $f:X\to\p^n$ be a regular, integral, divisorial double covering factoring through $P_m$ such that $Pic(X)\cong \Z f^*\osn(1)$ and let $\e$ be a rank $2$, $f$-aCM initialised vector bundle.
Then, either $\e\cong \os_X\oplus\os_X(\gamma)$with $\gamma\le0$ or any global section of $\e$ gives an extension as in
\sesl{\os_X}{\e}{\id_Y(\gamma)}{ext0}
for some $Y\subset X$ defined over $\gk$ such that $f|_Y$ an isomorphism onto a complete intersection $Z$ of type $(\alpha,\beta)$ for some integers $\alpha,\beta,\gamma$ such that $1\le\alpha,\beta\le m$ and $\gamma=\alpha+\beta-m$.
\One
\Oss
The hypothesis on the Picard group is automatic if $n\ge 4$ and $X$ cyclic or $n\ge 3$ and $\gk=\C$ by \autoref{Piccover}.
Moreover, if $n=2$ and $\gk=\C$ it holds for a general $X$ whenever $m\ge 3$.
When $Pic(X)$ is not cyclic or when $X$ is not regular, we expect more aCM sheaves to appear in addition to the above ones.
For example, as shown in \cite{Orlovmatrix}, the corresponding matrix factorisations determine the singularity category of $X$.
\one
\begin{proof}
Consider some section of $\e$, which exists since this sheaf is initialised.
Being $X$ regular, $Pic(X)\cong \Z f^*\osn(1)$ and $h^0(\e(-1))=0$, the vanishing locus of any section cannot be a divisor in $X$ and hence is empty or of codimension $2$.
In the first case, the cokernel is another line bundle, call it $\os_X(\gamma)$.
Being $\e$ initialised we have $\gamma\le0$.
Then $ext^1(\os_X(\gamma),\os_X))=h^1(\os_X(-\gamma))=0$ hence the extension is split.

In the second case, we must have $n\ge 2$.
Call $Y$ the codimension $2$ zero locus of some section of $\e$, so that we have the sequence \eqref{ext0} with $\os_X(\gamma)\cong det(\e)$.
Draw the first two rows of the diagram 
     \begin{equation}\label{pushideale}
         \begin{tikzcd}
          0 \ar[r] & \id_Z \ar[r] \ar[d] & \osn \ar[r] \ar[d] & \os_Z \ar[r] \ar[d] & 0 \\
          0 \ar[r] & f_*\id_Y \ar[r] \ar[d] & \osn\oplus\osn(-m) \ar[d] \ar[r]  & f_*\os_Y\ar[r] \ar[d] & 0 \\
         & \f \ar[r] & \osn(-m) \ar[r] & \g &
         \end{tikzcd}
     \end{equation}
     where $Z$ is defined to be the schematic image of $Y$, hence $\id_Z=ker(\osn\to f_*\os_Y)$.
     Everything commutes and we define $\f$ to be the cokernel of the induced map $\id_Z\to f_*\id_Y$ and $\g$ the cokernel of $\os_Z\to f_*\os_Y$.
     From snake lemma it follows that the bottom sequence is exact, in particular $\g\cong \os_{Z'}(-m)$ and $\f\cong \id_{Z'}(-m)$ for some subscheme $Z'$ of $\p^n$.
     Applying $f_*$ to \eqref{ext0}, recalling \autoref{divisorialcover} and \autoref{noether}, we get
     \sesl{\osn\oplus\osn(-m)}{\bigoplus_{i=0}^3\osn(\alpha_i)}{f_*\id_Y(\gamma)}{ext0'}
    which is short and exact being $f$ finite.
    Assume $\alpha_0\ge\dots\ge\alpha_3$.
    Being $\e$ initialised we must have $\alpha_0=0$ hence, being any non-zero endomorphism of $\osn$ an automorphism, we can cancel out this term and get another short exact sequence
\sesl{\osn(-m)}{\bigoplus_{i=1}^3\osn(\alpha_i)}{f_*\id_Y(\gamma)}{ext0''} 
    Moreover, the last part of \cite[Lem. 3.2.1]{HuyLeh} implies $\mu(\e)= 2(\mu(f_*\e)-\mu(f_*\os_X))$ which gives $\sum_{i=1}^3\alpha_i=2\gamma-2m$.
    Composing the right morphism in \eqref{ext0''} with the left column in \eqref{pushideale} twisted by $\osn(\gamma)$, we get a surjection \newline $ \phi:\bigoplus_{i=1}^3\osn(\alpha_i)\twoheadrightarrow \id_{Z'}(\gamma-m)$ hence $\alpha_3\le \gamma-m$.
    We will prove that it is always an equality and $Z'=\emptyset$.
    Indeed, if $\alpha_2>\gamma-m$ then $\alpha_1>\gamma-m$ hence $hom(\osn(\alpha_i),\id_{Z'}(\gamma-m))=0$ for $i=1,2$.
    A surjective morphism $\osn(\alpha_3)\to\id_{Z'}(\gamma-m)$  must be an isomorphism hence $Z'=\emptyset$ and $\alpha_3= \gamma-m$.
    In case $\alpha_2\le \gamma-m$ we have
    \begin{equation}\label{ineq}
        2\gamma-2m=\sum_{i=1}^3\alpha_i\le\alpha_1-2\gamma-2m\le 2\gamma-2m.
    \end{equation}
    hence $\alpha_1=0$ and $\alpha_2= \gamma-m=\alpha_3$. 
    But then $Z'=\emptyset$ since otherwise $hom(\osn(\gamma-m),\id_{Z'}(\gamma-m))=0$.
    Therefore, $\os_Z\cong f_*\os_Y$ and $\f\cong \osn(-m)$.
    In particular, being $f|_Y$ finite, since $f$ such, we deduce that $f|_Y$ is an isomorphism on the image, which is exactly $Z$.
    Thus, we showed that $\phi$ is split hence also the left column of \eqref{pushideale} is split, giving $f_*\id_Y(\gamma)\cong \id_Z(\gamma)\oplus\osn(\alpha_3)$.
    Simplifying \eqref{ext0''} we deduce a short exact sequence
    \ses{\osn(-m)}{\osn(\alpha_1)\oplus\osn(\alpha_2)}{\id_Z(\gamma)}
    Define $\alpha:=\gamma-\alpha_1$ and $\beta:=\gamma-\alpha_2$.
    Then the above sequence twisted by $\osn(-\gamma)$ gives
    \ses{\osn(-m-\gamma)}{\osn(-\alpha)\oplus\osn(-\beta)}{\id_Z}
    in particular $\alpha+\beta=m+\gamma$. Therefore, $Z$ is an $(\alpha,\beta)$ complete intersection and $\alpha_i\ge\gamma-m$ imply $\alpha,\beta\le m$.
\end{proof}

In particular, in the non-split case, those bundles cannot be pullbacks of bundles on $\p^3$.
What we really seek is the converse to the above result, i.e. a way to construct those sheaves starting from special subvarieties.
We begin with some standard computations that will be useful later.

\Le\label{numericsci}
Let $Z\subset \p^n$ with $n\geq 2$ be an $(\alpha,\beta)$ complete intersection, then
\begin{itemize}
    \item $\omega_Z\cong\os_Z(\alpha+\beta-n-1)$
    \item $h^0(\id_Z(i))= h^0(\osn(i-\alpha))+ h^0(\osn(i-\beta))-h^0(\osn(i-\alpha-\beta))$
    \item $h^j(\id_Z(i))=0,\;\; i\in \Z,\; 1\leq j\leq n-2$
    \item $h^{n-1}(\id_Z(i+\alpha+\beta-m))-h^{n}(\id_Z(i+\alpha+\beta-m))=h^n(\osn(i-m))-h^n(\osn(i+\alpha-m))-h^n(\osn(i+\beta-m))$
    \item $h^0(\os_Z(i))=\binom{i+3}{3}-h^0(\id_Z(i))$ for $i\ge 0$.
\end{itemize}
\ma
\begin{proof}
    The first claim follows by adjunction formula and the fact that $\omega_{\p^n}\cong\osn(-n-1)$.
    By assumptions on $Z$ we have an exact sequence \eqref{ic}.
    We have $h^j(\p^n,\osn(i))=0$ for $i\in \Z, 1\leq j\leq n-1$ hence, from the above sequence it immediately follows that $h^j(\p^n,\id_Z(i))=0$ for all $i\in \Z$ and $1\leq j\leq n-2$ and the formula in case $j>n-2$.
    Similarly, we can easily compute $h^0(\id_Z(i))$ for all $i\in\Z$.
    The last claim follows using the previous results and the cohomology of the adequate twists of
    \ses{\id_Z}{\osn}{\os_Z}
\end{proof}

We are ready to prove that any $Y$ as the ones constructed in \autoref{rk2-section} actually gives us an aCM sheaf of rank $2$, even without regularity assumptions.
\Te\label{section-rk2}
Suppose $n\ge3$ and fix a field $\gk$.
Let $f:X\to\p^n$ be an integral, divisorial double covering factoring through $P_m$.
Consider $Y\subset X$ closed of codimension $2$ and defined over $\gk$ such that $f|_Y$ is an isomorphism on the image and $Z:=f(Y)$ is an $(\alpha,\beta)$ complete intersection, with $1\le\alpha,\beta\le m$.
Any such $Y$ is the zero locus of a unique (up to scalar) section of a rank $2$, initialised, $f$-aCM sheaf $\e$ sitting in 
\sesl{\os_X}{\e}{\id_Y(\alpha+\beta-m)}{ext}
Moreover, $\e$ is reflexive, $\e\cong \e\du\otimes det(\e)$ and $(2m-\alpha-\beta)-n$ is the minimum $i\in\Z$ such that $h^n(\e(i))=0$.
\Ma
\begin{proof} 
We set $\gamma:=\alpha+\beta-m$, in particular $\gamma\le m$.
By \autoref{numericsci}, \autoref{divisorialcover} and $Z\cong Y$ we have that 
\[\omega_Y(-\gamma)\otimes\omega_X\du\cong\omega_Y(-\alpha-\beta+m)\otimes\omega_X\du\cong \os_Y(m-n-1)\otimes\omega_X\du\cong \os_Y.\]
Being $n\ge3$ we have $dim(Z)\ge 1$ and by \autoref{numericsci} $h^0(Y,\os_Y)=h^0(Z,\os_Z)=1$, since $Z$ is connected.
By \autoref{divisorialcover} $3)$ we know that $(X,\os_X(1))$ has no intermediate cohomology, in particular $h^i(X,\os_X(-\gamma))=0$ for $i=1,2$ since $n\ge 3$.
Therefore, applying $Hom(-,\os_X)$ to 
\ses{\id_Y(\gamma)}{\os_X(\gamma)}{\os_Z(\gamma)}
and using Serre duality twice we get
\[ext^1(\id_Y(\gamma),\os_X)\cong ext^2(\os_Y(\gamma),\os_X)\cong H^{n-2}(\os_Y(\gamma)\otimes\omega_X))\cong h^0(\omega_Y(-\gamma)\otimes\omega_X))\cong h^0(\os_Y)=1,\]
that is, we have a unique (up to isomorphism) non-split extension like \eqref{ext}.
Being $f_*\os_X\cong \osn\oplus\osn(-m)$, applying $f_*$ to \eqref{ext} we get
\sesl{\osn\oplus\osn(-m)}{f_*\e}{f_*\id_Y(\gamma)}{ext'}
We can draw the diagram \eqref{pushideale} and now, being $Y\cong Z$ we have $\g=0$ hence $\f=\osn(-m)$, the left column reads
\sesl{\id_Z(\gamma)}{f_*\id_Y(\gamma)}{\osn(\gamma-m)}{pid}
By \eqref{pid} and \autoref{numericsci} we have $h^j(f_*\id_Y(i))=0$ for $i\in\Z$ if $1\le j\le n-2$ and $h^0(f_*\id_Y(\gamma-1))=0$ hence, due to the fact that $f$ is finite and hence preserves cohomology, from \eqref{ext'} we get $h^j(X,\e(i))=0$ for $i\in\Z$ if $1\le j\le n-2$ and $\e$ initialised.
To prove that $\e$ is aCM we are left to show $h^{n-1}(X,\e(i))=h^{n-1}(f_*\e(i))=0$ for all $i\in\Z$.
By \eqref{ext'}, \eqref{pid} and \autoref{numericsci} we have
\[h^{n-1}(\e(i))=h^{n-1}(\id_Z(\gamma+i)-h^{n}(\osn(i))-h^{n}(\osn(i-m))+h^n(\e(i))-h^{n}(\id_Z(\gamma+i))-h^{n}(\osn(i+\gamma-m)).\]
Note that $h^{n-1}(\id_Z(i+\gamma))-h^{n}(\id_Z(i+\gamma))$ has already been computed in \autoref{numericsci} while to compute $h^n(\e(i))$ will require more work.
Apply $\Hom(-,\os_X)$ to \eqref{ext}: we get an exact sequence which can be split into
\[0\to\os_X(-\gamma)\to\Hom(\e,\os_X)\to\ck\to 0 \qquad \qquad 0\to \ck\to\os_X\to\ext^1(\id_Y(\gamma),\os_X)\to 0\]
which are everywhere exact except, possibly, for the right side of the second.
Since the extension class of \eqref{ext} is non-zero by definition, it follows that the aforementioned map is non-zero and hence surjective, being
\[\ext^1(\id_Y(\gamma),\os_X)\cong \ext^2(\os_Y(\gamma),\os_X)\cong\omega_Y(-\gamma)\otimes\omega_X\du\cong \os_Y\]
and $Y$ connected.
It follows that actually $\ck\cong \id_Y$ and hence the first sequence is a non-split extension of $\id_Y$ and $\os_X(-\gamma)$, i.e. $\e\du\cong \e(-\gamma)$ by the unicity proved before.
We deduce immediately that $\e\dd\cong (\e(-\gamma))\du\cong \e$, hence $\e$ is reflexive.
Moreover, by Serre duality we have
\[h^n(\e(i))=hom(\e(i),\omega_X)=h^0(\e\du(m-n-1-i))=h^0(\e(m-n-1-i-\gamma))\]
which, being $\e$ initialised, vanishes if and only if $i\ge m-n-\gamma=2m-\alpha-\beta-n$.
Recalling $h^0(\e(i))=h^0(f_*\e(i))$, by \eqref{ext'} and Serre duality the last equality becomes
\[h^n(\e(i))=h^n(\osn(\gamma+i-m))+h^n(\osn(\gamma+i))+h^n(\osn(i))+h^0(\id_Z(m-n-1-i)).\]
Finally, \autoref{numericsci} and Serre duality imply
\[h^0(\id_Z(m-n-1-i))=h^n(\osn(i+\alpha-m))+h^n(\osn(i+\beta-m))-h^n(\osn(i+\gamma))\]
therefore, putting everything together in the formula for $h^{n-1}(\e(i))$ we conclude that $\e$ is aCM.
\end{proof}

As a consequence, we obtain a general description of the splitting type of $\e$ only depending on $\alpha,\beta$, in particular this allows us to distinguish Ulrich bundles among other aCM-s.
\Co\label{splittype}
We have $f_*\e\cong \osn\oplus\osn(\alpha-m)\oplus\osn(\beta-m)\oplus\osn(\alpha+\beta-2m)$, in particular $\e$ is Ulrich if and only if $\alpha=\beta=m$.
Moreover, we have
\[ h^0(\e)=\begin{cases} 1  \hfill \alpha\le\beta<m \\
2 \qquad \hfill   \alpha<\beta=m \\
4 \hfill   \alpha=\beta=m \\
\end{cases}\]
\io
\begin{proof}
    We know that $f_*\e$ is split by \autoref{noether}.
    All the summands are of the form $\osn(\alpha_i)$ and, being $\e$ initialised, we have $\alpha_i\le 0$ and there is at least one $i$ such that $\alpha_i=0$.
    The other ones can be detected understanding when the function $h^0(\e(i))$ jumps.
    From \eqref{ext'} we have
    \[h^0(\e(i))=h^0(\osn(i))+h^0(\osn(i-m))+h^0(\osn(i+\gamma-m))+h^0(\id_Z(i+\gamma))=\]
    \[=h^0(\osn(i))+h^0(\osn(i+\alpha+\beta-2m))+h^0(\osn(i+\beta-m))+h^0(\osn(i+\alpha-m))\]
    where in the last step we used \autoref{numericsci}.
    From \autoref{noether} $\e$ is Ulrich if and only if all the summands of $f_*\e$ are trivial, hence $\alpha=\beta=m$.
    Finally, putting $i=0$ we get the last claim.
\end{proof}

\Co\label{condstab}
Suppose $X$ is regular and $Pic(X)\cong \Z H$.
Then the vector bundles constructed in \autoref{section-rk2} are 
\begin{itemize}
    \item slope-(semi)stable if and only if $\alpha+\beta>(\ge) m$
    \item not Gieseker-semistable if $\alpha+\beta=m$
    \item not simple if $\alpha+\beta\le m$.
\end{itemize}.
\io
\begin{proof}
If $\e$ is slope-(semi)stable then by \eqref{ext} we immediately get 
\[\dfrac{\alpha+\beta-m}{2}=\mu(\e)>(\ge) \mu(\os_X)=0.\]
For the converse, note that $\e$ as in \eqref{ext} is locally free, being $X$ regular.
Assume by contradiction that there is a slope-destabilising subsheaf $\f\subset\e$.
Then $\f$ must be torsion-free and of rank $1$ but, being $\e$ locally free, we can assume, up to taking its double dual, that $\f$ is reflexive and hence a line bundle, since $X$ is regular.
Now, our assumption on $Pic(X)$ implies $\f\cong\os_X(i)$, so that $0\ne hom(\os_X(i),\e)\cong h^0(\e(-i))$ and then $i\le 0$, being $\e$ initialised, which contradicts $\alpha+\beta>(\ge) m$ by the above slope computation.
Note that for $\alpha+\beta=m$ the sequence \eqref{ext} is Gieseker-destabilising since the Hilbert polynomial of $\os_X$ is bigger then the one of $\id_Y$.
Finally, if $\alpha+\beta\le m$ we get a non-trivial endomorphism of $\e$ as follows: compose the surjection $\e\twoheadrightarrow\id_Y(\gamma)$ given by \eqref{ext} with the inclusion $\id_Y(\gamma)\hookrightarrow \os_X(\gamma)\hookrightarrow\e$.
\end{proof}

We now put together the material of this section.
Note that, the claim $(1)\Rightarrow (2)$ in \autoref{2} is a consequence of \autoref{rk2-section} while $(2)\Rightarrow (1)$ follows from \autoref{section-rk2}.
The fact that the case $\alpha=m=\beta$ corresponds to Ulrich bundles follows from \autoref{splittype}.
The next lemma proves $(2)\Leftrightarrow (3)$.
Actually, the same ideas can be used quite generally to characterise subvarieties $Y\subset X$ of double coverings $f:X\to \p^n$ mapped isomorphically by $f$, see \cite[§3.1.5]{tesi}

\Le\label{23}
Supposing $char(\gk)\ne 2$ and $n\ge3$, the implications $(2)\Leftrightarrow (3)$ in \autoref{2} hold for any integral, cyclic double covering.
\ma
\begin{proof}
Let us start with $(2)\Rightarrow (3)$.
We need to prove that, if there is some $Y\subset X$ mapped isomorphically by $f$ to a complete intersection $Z$ of type $(\alpha,\beta)$ then we can write $b=p_\alpha q_{\alpha}+p_{\beta} q_{\beta}+p_m^2$ for $p_l\in H^0(\p^n,\osn(l))$, where $b=0$ defines the branch locus $B\subset\p^n$ of $X$.
Call $p_\alpha,p_\beta$ the two polynomials of degree $\alpha,\beta$ defining $Z$.
If $Z\subset B$ then we can write $b=p_\alpha q_{\alpha}+p_{\beta} q_{\beta}$ for some homogeneous polynomials $q_\alpha,q_\beta$ of degrees $2m-\alpha,2m-\beta$ respectively and we are done.
Otherwise, call $i_Z:Z\hookrightarrow \p^n$ and $i_Y:Y\hookrightarrow Y$, so that $f\circ i_Y=i_Z$ identifying $Z$ and $Y$.
Call $R$ the ramification locus of $f$, so that $f^*B=2R$.
Functoriality of pullback implies that 
    \[B|_Z=i_Z^*B=i_Y^*(f^*B)=i_Y^*(2R)=2i_Y^*R.\]
    Finally, notice that
    \[\os_Z(i_Y^*R)\cong i_Y^*\os_X(m)\cong i_Y^*f^*\osn(m)\cong\os_Z(m)\]
    hence $B|_Z$ is twice a divisor in the linear system $|\os_Z(m)|$ which, since complete intersections are projectively normal, is cut by an hypersurface of degree $m$ in $\p^n$.
    In other words, there exists some homogeneous polynomial $p_m$ of degree $m$ such that $b|_Z=p_m^2|_Z$.
    Therefore, $b=p_\alpha q_{\alpha}+p_{\beta} q_{\beta}+p_m^2$ for some homogeneous polynomials $q_\alpha,q_\beta$ of degrees $2m-\alpha,2m-\beta$ respectively.

    We are left to prove the converse.
    Recall that $X$ sits inside $P_m=\p(\osn\oplus\osn(-m))$ with equation $t^2-b=0$ and suppose that $b=p_\alpha q_{\alpha}+p_{\beta} q_{\beta}-p_m^2$.
    Define $Z=\{p_\alpha=0=p_\beta\}$.
    If we call $P_m\times_{\p^n}Z=\p(\os_Z\oplus\os_Z(-m))=:P_Z$ then the scheme $f\inv(Z)$ has equation $t^2-b=t^2-p_m^2=(t-p_m)(t+p_m)$ in $P_Z$.
    It follows that $f\inv(Z)$ has two components corresponding to $t\pm p_m=0$; we work out the case $t+p_m$ and call $Y$ the corresponding component.
    This polynomial defines a morphism
    \[0\to \os_{P_Z}\xrightarrow{t+p_m}\os_{P_Z}(1)\to \os_{P_Z}(1)\otimes \os_{Y}\to 0\]
    whose pushforward, using \cite[Prop. 22.86]{GW2}, is
    \[0\to \os_{Z}\xrightarrow{(1,p_m)}\os_{Z}\oplus\os_Z(m)\to \os_Z(m)\to 0\]
    and $Y$ is the section of $\pi|_{P_Z}$ over $Z$ corresponding to the surjection on the right.
    It follows that $\pi|_{Y}=f|_{Y}:Y\to Z$ is an isomorphism.
\end{proof}

For future use, we will show that $Y$ mapped isomorphically by $f$ to $Z$ comes exactly from the construction described in the second part of the above proof.

\Le\label{icp}
Let $f:X\to\p^n$ be an integral, cyclic double covering and $char(\gk)\ne2$. 
Suppose $Y\subset X\subset P_m$ is mapped isomorphically by $f$ to a complete intersection $Z$ of type $\alpha,\beta$, then $Y$ is a complete intersection in $P_m$ and $\n_{Y/P_m}\cong \os_Y(\alpha)\oplus\os_Y(\beta)\oplus\os_Y(m)$.
\ma
\begin{proof}
    Suppose that $Z\subset\p^3$ is the complete intersection of two divisors $D_\alpha,D_\beta$. We can view $Y$ as a section of $\pi_Z:\p_Z\to Z$, where $P_Z:=\p(\os_Z\oplus\os_Z(m))$ is naturally viewed as a subscheme of $P_m$.
If we call $\iota$ the covering involution of $X$, then also the variety $\iota(Y)$ must be a section of $\pi_Z$ and $X\cap P_Z=Y\cup \iota(Y)$ since $deg(f)=2$.
Moreover, if the branch locus of $f$ has equation $b=0$, then the equation of $X$ in $P_m$ is $t^2-b=0$, with $t$ the tautological section of $\os_{P_m}(1)$.
Thus, the only possibility of $X|_{P_Z}$ being reducible is that $t^2-b$ factors when restricted to $P_Z$.
It follows that $Y$ and $\iota(Y)$ have equation $t\pm\mathfrak{b}$, where $\mathfrak{b}^2=b|_{P_Z}$, in particular both are divisors in $|\os_{P_Z}(1)|$ mapping isomorphically to $Z$.
Note that the surjection $\psi:\os_{P_m}(1)\to\os_{P_Z}(1)$ induces a surjective map $H^0(\os_{P_m}(1))\to H^0(\os_{P_Z}(1))$: indeed, applying $\pi_*$ to $\psi$ we get the standard projection $\on{3}\oplus\on{3}(m)\to \os_Z\oplus\os_Z(m)$ which is surjective on global section being $Z$ a complete intersection, see \autoref{numericsci}.
We deduce that $Y$ is a complete intersection in $P_m$ of $\pi^*D_\alpha,\pi^*D_\beta$ and some divisor in $|\os_{P_m}(1)|$.
By \autoref{divisorialcover}, we get $\os_{P_m}(1)|_Y\cong \os_Y(m)$.
\end{proof}

\subsection{Relation with involution}

Since a cyclic double covering has an involution $\iota$, we want to investigate its behaviour on aCM sheaves.
Note that, if $\e$ is a direct sum of line bundles of the form $\os_X(\gamma)$, then it is fixed by $\iota$ since we have $\iota^*\os_X(1)\cong\os_X(1)$ being $\os_X(1)=f^*\osn(1)$.

\Prop\label{fisso}
Suppose $char(\gk)\neq 2$ and let $f:X\to \p^n$ be an integral, cyclic double covering in $P_m$ with $n\geq 3$.
Let $\e$ be a rank $2$, non-split initialised $f$-aCM bundle on $X$ and $\iota$ the covering involution.
Consider the statements:
\begin{enumerate}[i)]
    \item we have an exact sequence
    \sesl{\os_X}{\e}{\id_Y(\gamma)}{ext''}
    with $Y=\iota^* Y$ mapping isomorphically to a complete intersection of type $\alpha,\beta$  in $\p^n$ where $\gamma=\alpha+\beta-m$
    \item $\e\cong \iota^*\e$
\end{enumerate}
We have $i)\Rightarrow ii)$ and, if we also assume $X$ regular and $Pic(X)\cong \Z \os_X(1)$, then $ii)\Rightarrow i)$.
\One
\begin{proof}
\gel{$\mathbf{i)\Rightarrow ii}$}  
$Y=\iota^* Y$ implies $\iota^*\id_Y\cong \id_{\iota^*Y}\cong \id_Y$ hence applying $\iota^*$ to $\eqref{ext''}$  we get 
\sesl{\os_X}{\iota^*\e}{\id_Y(\gamma)}{ext'''}
Being $\e$ and $\iota^*\e$ aCM, the sequences \eqref{ext''} and \eqref{ext'''} are non-split but, as shown in \autoref{section-rk2}, we have $ext^1(\id_Y(\gamma),\os_X)=1$ hence the desired isomorphism follows.

\gel{$\mathbf{ii)\Rightarrow i}$}  
   If $\e\cong \iota^*\e$, the automorphism $\iota$ acts linearly on $H^0(X,\e)$.
   Recall that, being $\iota$ an involution and $char(\gk)\neq 2$, the eigenvalues can only be $\pm 1$ and it is diagonalisable.
   Therefore, at least one of the eigenspaces is non-empty.
   Let us choose an eigensection $\os_X\to\e$.
   Due to our enhanced assumptions, we can apply \autoref{rk2-section} and get the sequence \eqref{ext''}.
Moreover, by construction, $\iota^*$ preserves the left morphism in this sequence, hence also its cokernel. 
We just proved that $\iota^*$ induces an isomorphism $\id_Y(\gamma)\cong \iota^*\id_Y(\gamma)\cong \id_{\iota^*Y}(\gamma)$ which implies $\iota^*Y=Y$.
\end{proof}

\Co\label{fisso'}
With the above notation, suppose $X$ regular, $Pic(X)$ is generated by $\os_X(1)$, $char(\gk)\neq 2$ and $\gk$ algebraically closed.
Then, there exists an $f$-aCM, non-split, initialised bundle $\e$ of rank $2$ such that $\e\cong \iota^*\e$ if and only if there is some $Y$ contained in the ramification divisor of $f$ such that \eqref{ext''} holds.
\io
\begin{proof}
    We can apply both conclusions in \autoref{fisso} hence, $\e\cong \iota^*\e$ if and only if there exists $Y\subset X$ of codimension $2$ such that $Y=\iota^* Y$ is the zero locus of a section of $\e$, as in \eqref{ext''}.
    But from \autoref{rk2-section} $f|_Y$ is injective hence $Y$ must be pointwise fixed by $\iota$, that is $Y\subset R$.
    Conversely, if $Y\subset R$ then clearly $Y=\iota^* Y$.
\end{proof}

Actually, Ulrich sheaf are globally generated, \autoref{corollario}, and their involutes naturally appear as syzygy sheaves, up to twist.
\Prop\label{Discount}
Suppose $char(\gk)\neq 2$.
Let $f:X\to\p^n$ be an integral, cyclic double covering with branch locus of degree $2m$ and call $\iota$ the involution of this covering.
If $\e$ is an $f$-Ulrich sheaf of rank $r$ then, taking the evaluation morphism of global sections of $\e$ we get
\sesl{\f(-m)}{\os_X^{2r}}{\e}{discount}
where $\f\cong \iota^*\e$.
\One
\begin{proof}
Recall that $X$ can be seen as a divisor of equation $p=t^2-b$ inside $P_m$ for some $b\in H^0(\p^n,\osn(2m))$, in particular $X\in |\os_{P_m}(2)|$.
In this case, the matrix $A$ appearing in \autoref{factorulrich} has coefficients in $H^0(P_m,\os_{P_m}(1))$ hence, up to changing of bases we can assume it is of the form $t\gi+A'$ where the coefficients of $A'$ are all in $H^0(\p^n,\osn(m))$.
Adding \eqref{discount} we can form the commutative diagram below
\dia
0 \ar[r] & \os_{P_m}(-2)^{2r} \ar[r,"id"] \ar[d,"M"] & \os_{P_m}(-2)^{2r} \ar[d,"p\gi_{2r}"] & &  \\
    0 \ar[r] & \os_{P_m}(-1)^{2r} \ar[d] \ar[r,"A=t\gi_{2r}+A'"] &  \os_{P_m}^{2r} \ar[d] \ar[r] & \e  \ar[r] \ar[d,"id"] & 0 \\
    0 \ar[r] & \f(-m) \ar[r]  &  \os_X^{2r} \ar[r] & \e  \ar[r] & 0 \\
    \mma
where $M$ is defined to be the matrix with coefficients in $H^0(P,\os_{P_m}(1))$ resulting from the diagram.
Twisting the left column by $\os_{P_m}(1)$, and recalling that $\os_{P_m}(1)|_X\cong \os_X(m)$ by \autoref{divisorialcover}, we get immediately that $\f$ is Ulrich by \autoref{Ulrichdivisorial}.
We can write $M=tM_1+M_2$ where $M_i$ have coefficients in $H^0(\p^n,\osn(m))$, in particular do not contain $t$.
Commutativity of the higher square implies that
\[(t^2-b)\gi_{2r}=p\gi_{2r}\circ id_{\os_P(-2)^{2r}}=(t\gi_{2r}+A')(tM_1+M_2)=t^2M_1+t(M_2+A'M_1)+A'M_2\]
therefore we get $M_1=\gi_{2r}$ and $M_2=-A'$ hence $M=tM_1-M_2$.
Recall that, $\iota$ is exactly the restriction to $X$ of the automorphism $\iota_P$ of $P$ given by sending $t\mt-t$.
But then
\[\iota_P(A)=\iota_P(t\gi_{2r}+A')=-t\gi_{2r}+A'=-M.\]
and since $M$ and $-M$ differ for just a sign, the corresponding cokernels are isomorphic hence $\f\cong\iota^*\e$.
\end{proof}

\subsection{Positivity properties}\label{+}

For simplicity, in this subsection we assume $\gk=\C$.
Here we want to understand the positivity properties of a special $f$-Ulrich bundle $\e$ of rank $2$ on a smooth double covering $f:X\to\p^n$ with $n\ge 3$ and $Pic(X)\cong \Z H$, where $H=c_1(\os_X(1))$.
We already know that $\e$ is globally generated by \autoref{corollario}.
Moreover, notice that by the computation of $h^0(\e)$ given in \autoref{splittype}, among the aCM bundles described in \autoref{section-rk2} only the Ulrich ones are globally generated.
We will show that the codimension $2$ subschemes in $X$ on which $\e$ restricts to a non-ample bundle are exactly the zero loci of sections of $\iota^*\e$, where $\iota$ is the involution on $X$ associated to $f$.

Consider the projective bundle $q:\p(\e)\to X$.
Define $\xi:=c_1(\os_{\p(\e)}(1))$ and $L:=c_1(q^*\os_X(1))$.
Being $Pic(X)\cong \Z H$, we have $Pic(\p(\e))=<\xi,L>$ by \cite[Prop. 24.69]{GW2}, moreover linear and numerical equivalence coincide on those varieties.
Recall that a divisor $D$ is \textbf{Nef} if $D\cdot C=deg(D|_C)\ge 0$ for all curves $C\subset X$; they form a convex cone in $Pic(X)\otimes\Q$.
Also effective divisors on $X$ form a cone and its closure in $Pic(X)\otimes\Q$ is called the \textbf{Pseudo-effective} cone.
Divisors which are in the interior of this cone are called \textbf{big}. 
We can easily describe both cones for $\p(\e)$.
\Le\label{coni}
In the above setting, for any rank $2$ Ulrich bundle $\e$ on $X$ we have $\xi^4=0$, in particular $\e$ is not big\footnote{Recall that a vector bundle $\e$ is called big if $\os_{\p(\e)}(1)$ is a big line bundle}.
Both Nef and pseudo-effective cone of divisors of $\p(\e)$ are spanned by $L,\xi$.
\ma
\begin{proof}
Since $\os_X(1),\e$ are globally generated, the divisors $L,\xi$ are basepoint-free.
It follows that they are both Nef.
By \cite[Thm. 2.2.16]{Laz1}, they are not big if we verify that their top self-intersection is $0$.
Clearly 
\[L^{n+1}\sim (q^*H)^{n+1}\sim q^*(H^{n+1})\sim 0.\]
Since $\xi$ is globally generated and $h^0(\p(\e),\os_{\p(\e)}(1))=h^0(X,\e)=4$ by \autoref{splittype}, its linear system induces a morphism $\phi:\p(\e)\to\p^3$ such that $\phi^*\on{3}(1)\cong \os_{\p(\e)}(1)$.
It follows that $\xi^4\sim 0$.
This implies that both $L,\xi$ are Nef not big hence are extremal in both the Nef and the pseudo-effective cone of $X$.
\end{proof}

As promised at the beginning, we come back to non-ample loci for the restriction of $\e$.
We show that those subvarieties lift to $\p(\e)$ and are contained in the fibers of the morphism determined by $\os_{\p(\e)}(1)$.
In order to do so, we will need to understand the restriction of $\e$ to the zero loci of $\iota^*\e$.
We start with the following lemma.
\Le\label{intersezionebranch}
Let $f:X\to \p^n$ be a double covering with ramification divisor $R$. 
If $Z\subset\p^n$ is a variety such that $f\inv(Z)$ has two components $Z_i$ mapped isomorphically on $Z$ by $f$ then 
\[Z_1\cap R=Z_1\cap Z_2=Z_2\cap R.\]
\ma
\begin{proof}
Suppose $X$ has equation $t^2-b=0$ in $P_m$ and call $\pi:P_m\to \p^n$ the standard projection.
By assumption $\pi\inv(Z)$ cuts on $X$ exactly $Z_1\cup Z_2$, which means $Z_i$ have equations $t\pm\mathfrak{b}=0$ in $\pi\inv(Z)$, for some $\mathfrak{b}^2=b$.
But cutting $X$ with the hyperplane given by $t=0$ we obtain exactly the ramification divisor $R$.
In particular, the ideal in $X$ defining $Z_1\cap Z_2$ is the same as the one defining $Z_1\cap R$, hence the conclusion follows.
\end{proof}

\Le\label{restrizioneinv}
Suppose $\e$ is a rank $2$ Ulrich bundle on a smooth double covering of $\p^n$.
If $Y$ is a zero locus of some section of $\e$ then $\e|_{\iota^*Y}\cong \os_{\iota^*Y}\oplus\os_{\iota^*Y}(m)$.
\ma
\begin{proof}
    By \autoref{rk2-section}, the zero locus of a section of $\e$ is a codimension $2$ subscheme $Y$.
    By \autoref{fisso'} a general $Y$ is distinct from $\iota^*Y$, which is the zero locus of a section of $\iota^*\e$.
 By assumption we have
    \sesl{\os_X}{\e}{\id_Y(m)}{solita}
    Restricting 
    \ses{\id_Y(m)}{\os_X(m)}{\os_Y(m)}
    to $\iota^*Y$ we get the right exact sequence
    \ses{\id_Y(m)|_{\iota^*Y}}{\os_{\iota^*Y}(m)}{\os_{Y\cap \iota^*Y}(m)}
    hence we have a surjection from $\id_Y(m)|_{\iota^*Y}$ to 
    \[Ker\left(\os_{\iota^*Y}(m)\twoheadrightarrow\os_{Y\cap \iota^*Y}(m)\right)=\id_{Y\cap \iota^*Y/\iota^*Y}(m)\cong \os_{\iota^*Y}(-R)(m)\cong \os_{\iota^*Y}\]
    since we have $Y\cap \iota^*Y=R\cap\iota^*Y$, by \autoref{intersezionebranch}, and $\os_X(R)\cong\os_X(m)$.
    We derive a surjection \[\e\twoheadrightarrow\e|_{\iota^*Y}\twoheadrightarrow\id_Y(m)|_{\iota^*Y}\twoheadrightarrow\os_{\iota^*Y}.\]
    The kernel of the morphism $\e|_{\iota^*Y}\twoheadrightarrow \os_{\iota^*Y}$ has to be locally free of rank $1$ hence isomorphic to $det(\e|_{\iota^*Y})\cong\os_{\iota^*Y}(m)$, so we deduce the sequence
    \sesl{\os_{\iota^*Y}(m)}{\e|_{\iota^*Y}}{\os_{\iota^*Y}}{zeri}
    Since $\e$ is globally generated also $\e|_{\iota^*Y}$ must be such hence this sequence splits. 
\end{proof}

Finally, we are ready for our result.
\Prop\label{nonampio}
For $n\ge 3$, let $f:X\to\p^n$ be a smooth, complex double covering, $\e$ a $f$-Ulrich bundle of rank $2$ such that $\iota^*\e\neq \e$.
Call $q:\p(\e)\to X$ the standard projection and $\phi:\p(\e)\to\p^3$ the morphism determined by $|\os_{\p(\e)}(1)|$.
Then $\phi$ is surjective and the family $\phi:\p(\e)\to\p^3$ is isomorphic to the family $\psi:\mathcal{U}\to\p(H^0(X,\iota^*\e)\du)\cong \p^3$ of zero loci of sections of $\iota^*\e$, in particular $\phi$ is flat.
\One
\begin{proof}
\textbf{Part} $\mathbf{1}:$ \textbf{Surjectivity of} $\mathbf{\phi}$

    Before continuing, we will need some more notation.
    Remember that by \autoref{section-rk2} we have $c_1(\e)\sim mH$ hence $q^*c_1(\e)\sim mL$.
Call $C=q^*Y$ and $\iota^*C:=q^*\iota^*Y$.
Let us recall that the Chow ring $A^{\bullet}(\p(\e))$ is isomorphic to $\Z[\xi,L]/(\xi^2-mL\cdot\xi+\gamma)$ by \cite[Example 8.3.4]{Ful}\footnote{our $\p(\e)$ is actually $P(\e\du)$ with Fulton's notation}.
Identifying a $0$-cycle with its degree we have 
\[H^{n-2}\cdot Y=m^2 \qquad H^n=2,\]
the first by \autoref{rk2-section}, and the second since $|H|$ induces the double covering $f$.
From those we also deduce 
\[\xi\cdot L^n=\xi\cdot q^*(H^n)=2 \qquad \xi\cdot C\cdot L^{n-2}=\xi\cdot q^*(c_2(\e)\cdot H^{n-2})=m^2.\]
We have
    \begin{equation}\label{csi3}
        \xi^3\sim \xi\cdot \xi^2=\xi (mL\cdot\xi-C)=mL(mL\xi-C)-\xi\cdot C=\xi(m^2L^2-C)-mLC\neq 0
    \end{equation}
    by the structure of $A^{\bullet}(\p(\e))$ or, for example, because we can compute
    \[\xi^3\cdot L^{n-2}=(\xi(m^2L^2-C)-mLC)\cdot L^{n-2}=m^2\xi\cdot L^n-\xi\cdot C\cdot L^{n-2}-mL^{n-1}\cdot C=\]
    \begin{equation}\label{contocubo} 
    =2m^2-m^2-m\cdot q^*(H^{n-1}\cdot Y)=m^2.
    \end{equation}
    Then $\phi$ is surjective, since otherwise its image would have dimension at most $2$ and hence $\xi^3=\phi^*c_1(\on{3}(1))^3=0$.

\textbf{Part} $\mathbf{2: \phi=\psi}$

The proof is quite technical but the idea is simple: since $\e|_{\iota^*Y}$ has a trivial quotient, by \autoref{restrizioneinv} then we have a lift $\iota^*Y\cong \p(\os_{\iota^*Y})\subset \p(\e)$ which is contained in some fiber of $\phi$.
Now, we repeat the previous construction in families, in order to get a morphism $\cu\to\p(\e)$ which we prove is actually an isomorphism.

    Consider the flat family $\psi:\cu\to\p^3\cong\p(H^0(X,\iota^*\e)\du)$ whose fibers are the zero loci $\psi_t:=\iota^*Y_t$ of sections of $\iota^*\e$, where $t\in\p^3$.
    Note that $\cu$ is irreducible of dimension $n+1$ since $\iota^*Y_t$ have codimension $2$ in $X$ and the general one is irreducible being connected, since isomorphic to a complete intersection in $\p^n$, and smooth, being $\e$ globally generated.
    Call $e:\cu\to X$ the evaluation morphism sending each curve to its realisation in $X$.
    \dia
    & \p(\e) \ar[rd, "q"]\ar[ld, "\phi"] & &  \cu \ar[rd, "\psi"]\ar[ld, "e"] \ar[ll, "g"]&  \\
    \p^3 & & X & & \p^3 \\
    \mma
    Then, $\psi_*e^*(\e(-m))\cong \cl$ is a line bundle by Grauert's theorem, see \cite[Thm. 23.140]{GW2}, since $h^0(\psi_t,(e^*\e(-m))_t)=h^0(\iota^*Y_t,\e|_{\iota^*Y_t}(-m))=1$ for all $t\in\p^3$.
    By projection formula we have 
    \[\psi_*(e^*\e(-m)\otimes\psi^*\cl\du)\cong \cl\otimes\cl\du\cong\on{3}\]
    hence $h^0(\cu,e^*\e(-m)\otimes\psi^*\cl\du)=1$.
    In other words, we have a morphism $e^*\os_X(m)\otimes\psi^*\cl\to e^*\e$ whose cokernel is a line bundle since restricted to all the fibers of $\psi$ coincides with $\os_{\psi_t}$ by \eqref{zeri}.
    A determinant computation implies we have an exact sequence
    \sesl{e^*\os_X(m)\otimes\psi^*\cl}{e^*\e}{\psi^*\cl\du}{magico}
    
    By universal property of $\p(\e)\to X$, this gives us a morphism $g:\cu\to\p(\e)$ over $X$ such that $g^*\os_{\p(\e)}(1)\cong \psi^*\cl\du$.
    Our next goal is to show that $g$ is birational.
    Since $g^*\phi^*\on{3}(1)\cong g^*\os_{\p(\e)}(1)\cong \psi^*\cl\du$, then the only curves that the morphism $\psi\circ g$ could contract are fibers of $\psi$.
    But $g|_{\psi_t}$ is an isomorphism on the image, since $e=q\circ g$ and $e|_{\psi_t}$ is an isomorphism on the image, by definition of $e$.
    Therefore, $g$ is finite and the fibers of $\phi\circ g$ are finite unions of $\psi_{t_i}$-s.
    To show that it is generically injective it is enough to show that the general fiber of $\phi$ contains just one $g(\psi_{t})$.
    But, if $s\neq t$ then we have $e(\psi_s)\neq e(\psi_t)$ hence $g(\psi_s)\neq g(\psi_t)$.
    Suppose, by contradiction, that the general fiber of $\phi$ contains both $g(\psi_t), g(\psi_s)$ with $g(\psi_t)\neq g(\psi_s)$.
    Being $\phi$ surjective, its general fiber has dimension $n-2$ and as a cycle is rationally equivalent to $\xi^3$.
    By dimensional reasons, $g(\psi_t)$ and $g(\psi_s)$ are both irreducible components of this fiber therefore, we have $\xi^3\sim g(\psi_t)+g(\psi_s)+\eta$ for some effective cycle $\eta$.
   But then by \eqref{contocubo}
    \[m^2=\xi^3\cdot L^{n-2}=(g(\psi_t)+g(\psi_s)+\eta)\cdot L^{n-2}\ge (g(\psi_t)+g(\psi_s))\cdot L^{n-2}=\]
    so by projection formula it becomes
    \[=g_*(\psi_t+\psi_s)\cdot L^{n-2}=(\iota^*Y_t+\iota^*Y_s)\cdot H^{n-2}=2m^2,\]
    a contradiction.
    Therefore, $g$ is a generically injective morphism between two varieties, which is birational, since we are over $\C$.
    Being $\p(\e)$ smooth, $g$ must be an isomorphism by Zariski's main theorem, \cite[Cor. 12.88]{GW1}. 
\end{proof}

\Oss
In the quadric case, i.e. $m=1$ and $n=3,4$, we have already seen that $\e$ is a spinor bundle and the $Y$-s are isomorphic to $\p^{n-2}$.
For $n=3$, this second contraction $\phi$ is again a projective bundle and in \cite[(3.4)]{SzurekWisniewski}, using different methods, it has been shown that it is the projectivization of a \emph{null-correlation bundle}.
\one

We end by an unrelated remark.
\Oss
By the universal property of Grassmannians, see \cite[Chapter 8 §4]{GW1}, any rank $2$ Ulrich bundle on a double covering $f:X\to\p^n$ gives us a morphism $\Phi_\e:X\to Gr(2,4)$ such that $\Phi_\e^*\mathcal{Q}\cong \e$, where $\mathcal{Q}$ is the universal quotient bundle on it.
Note that $Gr(2,4)$ is actually the smooth quadric in $\p^5$ and $Q$ is a spinor bundle, see \cite[Def. 1.3]{Ottaviani_spinor}, which is Ulrich respect to the restriction of $\on{5}(1)$. 
The map $\Phi_\e$ cannot be an embedding for $n\ge 3$ and $m\ne 1$.
This is clear for $n\ge 4$, while for $n=3$ follows from the fact that, by Lefschetz hyperplane theorem, see \cite[XII Cor. 3.7]{SGA2}, any divisor in this quadric has the Picard group generated by the restriction of $\on{5}(1)$, which is very ample, while $Pic(X)$ contains an ample but not very ample divisor.
Furthermore, for $n=3$ we have $det(\e)\cong \os_X(m)$ very ample.
This implies that the claims in the second part of \cite[Thm. 1]{LopezSierra_geom} are not equivalent if we only assume the polarisation to be ample and globally generated instead of very ample.
\one

\section{Existence of rank $2$ aCM bundles on some double coverings of $\p^3$}

In the previous section, through Hartshorne--Serre correspondence, we reduced the existence of aCM bundles on double coverings to the existence of subvarieties mapping isomorphically to complete intersections on projective space.
In this subsection, we show the existence of such subvarieties for a general double covering of $\p^3$ with branch locus of degree $2m=4,6,8$, thus proving the existence of the corresponding aCM bundles.
Note that the case $m=1$ corresponds to quadric hypersurfaces, on which aCM bundles are understood due to \cite{Knorrer}, hence, even though our methods are able to recover the cited results, we will always ignore it.
A dimensional computation shows that for $n=3$ and $m\ge 5$ or $n\ge 4$ and $m\ge 2$ the general double coverage of $\p^n$ does not admit such bundles, demonstrating that our result is somehow sharp.

\subsection{Existence of rank $2$ aCM sheaves}

This is the main result of this section.

\Prop\label{polinomi}
Suppose $char(\gk)\ne 2$ and $\gk$ infinite.
Fix $n=3$ and $m=2,3,4$.
For any $1\le \alpha, \beta\le m$ the general element $b\in H^0(\p^3,\on{3}(2m))$ can be written as $b=p_\alpha q_\alpha+p_\beta q_\beta +p_m^2$ with $p_l\in H^0(\p^3,\on{3}(l))$ and $q_l\in H^0(\p^3,\on{3}(2m-l))$ for $l=\alpha,\beta,m$.
\One
\begin{proof}
    Set $V_l:=H^0(\p^3,\osn(l))$ and $V_{\alpha,\beta}:=V_\alpha\times V_{2m-\alpha}\times V_\beta\times V_{2m-\beta}\times V_m$
Consider the algebraic map
\[\phi:V_{\alpha,\beta}\to V_{2m} \qquad (p_\alpha,q_\alpha,p_\beta,q_\beta, p_m)\mapsto p_\alpha q_\alpha+p_\beta q_\beta +p_m^2,\]
we want to show that this is dominant.
Its differential, seen as a linear map $V_{\alpha,\beta}\to V_{2m}$, in the point $(p_\alpha,q_\alpha,p_\beta,q_\beta, p_m)$ can be written as
\begin{equation}\label{differential}
    d_{(p_\alpha,\dots, p_m)}\phi: (p_\alpha',q_\alpha',p_\beta',q_\beta', p_m')\mt p_\alpha q_\alpha'+p_\beta q_\beta'+p_\alpha' q_\alpha+p_\beta' q_\beta+2p_mp_m'.
\end{equation}
If we show that such linear morphism is surjective for some choice of $p_\alpha,q_\alpha,p_\beta,q_\beta, p_m$, then $\phi$ would be smooth in that point, by \cite{GW2}[Thm. 18.74], and hence dominant.
Note that, since $char(\gk)\ne 2$, then we can ignore the coefficient $2$ appearing in the differential.
Therefore, its enough to find $5$ polynomials of degrees $\alpha,2m-\alpha,\beta,2m-\beta,m$ such that the ideal they generate contains all $V_{2m}$.
This can be done with the help of \textit{Macaulay2}, see \autoref{appendix}.
\end{proof}

\Co\label{rk2esistenza}
Suppose $char(\gk)\ne 2$, $\gk$ is an infinite field and $f:X\to\p^3$ is an integral, cyclic double covering branched along a surface of degree $2m=4,6,8$.
If $X$ is general enough then, for any $1\le\alpha,\beta\le 4$ it admits curves $Y\subset X$ mapping isomorphically to complete intersections of type $\alpha,\beta$ and hence the corresponding non-split rank $2$ aCM sheaves. 
\io
\begin{proof}
    Given \autoref{polinomi}, the first claim follows from \autoref{23} while the second from \autoref{ulrichsommeprodotti}.
\end{proof}

In most cases, we can refine this result by removing the "general" assumption.

\Prop\label{a=b}
If $\alpha\ne \beta$ or $\alpha=1=\beta$ then the conclusion of \autoref{rk2esistenza} is true for any integral $X$ as above.  
\One
\begin{proof}
    Consider the Hilbert scheme containing complete intersections $Z$ of type $\alpha,\beta$ in $\p^3$ and define $\mathcal{H_Z}$ to be the locus whose points actually correspond to complete intersections.
    Let us recall the explicit construction of this scheme given in \cite[§ 2.2.2]{Benoist}, which in particular shows that $\mathcal{H_Z}$ is a regular and, under our assumptions, projective variety.
Indeed, if $\alpha<\beta$, then $\mathcal{H_Z}$ is isomorphic to a projective bundle with base $|\on{3}(\alpha)|$ and fiber $|\os_W(\beta)|$ over the point $[W]\in |\on{3}(\alpha)|$.
Otherwise, if $\alpha=1=\beta$ then we have $\mathcal{H_Z}=Gr(2,H^0(\on{3}(1)))$.
On the ruled surface $\pi\inv(Z)$ we have the linear system $|\os_{P_m}(1)|_{\pi\inv(Z)}|$, which contains an open subset parametrising curves $Y$ mapping isomorphically to $Z$.
We consider the projective bundle $W\to\mathcal{H_Z}$ having such linear systems as fibers.
We define $\mathscr{I}$ as the projective bundle over $W$ having $|\id_{Y/P_m}(2)|$ as fiber over the point $Y$. 
In particular, $\mathscr{I}$ is a regular and projective variety.
Recall that $X$ can be seen as a divisor in $|\os_{P_m}(2)|$ and we call $\pi$ the projection $P_m\to\p^3$.
The variety $\mathscr{I}$ parametrises pairs $(Y,X)$ such that $\pi(Y)\in \mathcal{H_Z}$, $X\in |\os_{P_m}(2)|$ and $Y\subset X$, in particular we get the diagram
\[\xymatrix{ 
& \mathscr{I} \ar[ld]^{\pi_1} \ar[rd]^{\pi_2} & \\
\mathcal{H_Z} & & |\os_{P_m}(2)|}\]
But now $\pi_2$ is dominant by \autoref{rk2esistenza} and $\mathscr{I}$ is proper, so that $\pi_2$ is surjective.
Pick some integral $X\in |\os_X(2)|$ and some $Y\subset X$ such that $(Y,X)\in\mathscr{I}$.
Being $X$ integral, the map $\pi|_X:X\to\p^3$ is finite, then a fortiori $\pi|_Y:Y\to Z$ is finite.
Being $Y\in |\os_{P_m}(1)|_{\pi\inv(Z)}|$, we conclude that $\pi|_Y$ must be an isomorphism on the image.

\end{proof}

\Oss
As a consequence of the above proof, we have seen that $\mathscr{I}$ is regular.
If $char(\gk)=0$ then the general fiber of $\pi_2$ would be smooth.
If $\pi_2$ were smooth outside a codimension $2$ locus, then we could prove that this general fiber is also connected and hence irreducible.
Indeed, considering its Stein factorisation 
\[\mathscr{I}\xrightarrow{g} S\xrightarrow{h} |\os_{P_m}(2)|.\]
If we show that $h$ is an isomorphism then $\pi_2=g$ has connected fibers.
The morphism $h$ is finite but, being $\p^N$ simply connected, if it is not an isomorphism, then it has to ramify.
By purity of branch locus, \cite[\href{https://stacks.math.columbia.edu/tag/0BMB}{Tag 0BMB}]{Stacks}, the ramification locus must have codimension $1$, contradicting the smoothness of $\pi_2$ in codimension $2$.
\one

\subsection{Hilbert schemes and moduli spaces}

Our previous argument allows us to identify the component in the Hilbert scheme on $X$ that contains the curves $Y$ as the fibre of the morphism $\pi_2$ on $[X]$.
We will study its properties through the normal sheaf of $Y$ in $X$, using deformation theory.

\Le\label{cohonormale}
Suppose that $char(\gk)\ne 2$ and $\gk$ infinite.
Let $f:X\to \p^3$ be an integral, regular, cyclic double covering sending $Y\subset X$ isomorphically to some complete intersection $Z$ of type $\alpha,\beta$.
Then we have an exact sequence
\sesl{\n_{Y/X}}{\os_Y(m)\oplus\os_Y(\alpha)\oplus\os_Y(\beta)}{\os_Y(2m)}{normali}
Moreover, if $Y$ and $X$ are sufficiently general, $1\le \alpha\le \beta\le m$ and $m=2,3,4$ then
\[h^0(\n_{Y/X})= h^0(\os_Z(\alpha))+h^0(\os_Z(\beta))+h^0(\os_Z(m))-h^0(\os_Z(2m))=\alpha\cdot \beta\cdot (4-m)+h^1(\n_{Y/X})\]
\[ h^1(\n_{Y/X})=\begin{cases} 0  \hfill \alpha\le\beta<4 \\
1 \qquad \hfill   \alpha<\beta=m=4 \\
3 \hfill   \alpha=\beta=m=4 \\
\end{cases}\]
\ma
\begin{proof}
\textbf{Step 1: Setup}

\noindent Recall that $X$ can be seen as the zero locus of a section $p\in H^0(\os_{P_m}(2))$ of the $p=t^2-b$, where $P_m:=\p(\on{3}\oplus\on{3}(m))$.
Call $D_\alpha,D_\beta$ the two surfaces, given by polynomials $p_\alpha,p_\beta$, whose complete intersection is $Z$.
By \autoref{23} we have $b=p_\alpha q_{\alpha}+p_{\beta} q_{\beta}-p_m^2$, in particular we can write $p=(t+p_m)(t+p_m)-p_\alpha q_\alpha-p_\beta q_\beta$.
By \autoref{icp} $Y$ is a c.i. in $P_m$ of divisors defined by $p_\alpha,p_\beta, t+p_m$.
In particular, the conormal sheaf $\c_{Y/P_m}\cong \n_{Y/P_m}\du\cong \os_Y(-m)\oplus\os_Y(-\alpha)\oplus\os_Y(-\beta)$ is locally free.
The embeddings $Y\subset X\subset P_m$ give the conormal sequence 
\sesl{\c_{X/P_m}|_Y}{\c_{Y/P_m}}{\c_{Y/X}}{conormali}
which is exact by \cite[\href{https://stacks.math.columbia.edu/tag/06BA}{Tag 06BA}]{Stacks} being $Y\subset X$ l.c.i. by \cite[\href{https://stacks.math.columbia.edu/tag/0FJ2}{Tag 0FJ2}]{Stacks}.
Being all the sheaves in \eqref{conormali} locally free, we deduce that its dual sequence is still short exact and coincides with \eqref{normali}.
Taking cohomology we deduce $H^0(\n_{Y/X})= Ker(H^0(\n_{Y/P_m})\xrightarrow{\rho}H^0(\os_Y(2m))$, where we define $\rho$ to be the induced map on global sections.

\textbf{Step 2:} $\mathbf{\rho}$ \textbf{is surjective}

\noindent 
We can define a morphism $\os_{P_m}(-2)\to \os_{P_m}(-1)\oplus \pi^*\left(\on{3}(-\alpha)\oplus\on{3}(-\beta)\right)$ by the vector $(t-p_m,q_\alpha,q_\beta)$ and a morphism $\os_{P_m}(-1)\oplus \pi^*\left(\on{3}(-\alpha)\oplus\on{3}(-\beta)\right)\to \os_{P_m}$ by scalar product with $(t+p_m,-p_\alpha,-p_\beta)$.
Since the multiplication by $p$ gives $\os_{P_m}(-2)\cong \id_{X/P_m}\subset\os_{P_m}$, those maps fits in the diagram 
\dia
\os_{P_m}(-2) \ar[rrr, "\cdot p"'] \ar[d,"\ensuremath{(t-p_m,q_\alpha,q_\beta)}"'] & & & \id_{X/P_m} \ar[d,hookrightarrow] \\
\os_{P_m}(-1)\oplus \pi^*\left(\on{3}(-\alpha)\oplus\on{3}(-\beta)\right)\ar[rrr, twoheadrightarrow, "\ensuremath{\cdot (t+p_m,-p_\alpha,-p_\beta)}"] & & & \id_{Y/P_m} \\
\mma
which, once we apply $Hom(-,\os_Y)$, becomes 
\dia
 \os_Y(2m)   & & \n_{X/P_m}|_Y \ar[ll,"\cong"]  \\ \os_{Y}(m)\oplus \os_Y(\alpha)\oplus\os_Y(\beta) \ar[u,"\ensuremath{\cdot(t-p_m,q_\alpha,q_\beta)}"]  & & \n_{Y/P_m} \ar[u]  \ar[ll,"\cong"]  \\
\mma
Being $t+p_m=0$ on $Y$ we have $t-p_m=-2p_m$ on $Y$. Since $char(\gk)\ne 2$, for our aims we can replace $2p_m$ by $p_m$, hence the vertical maps are scalar product with $(-p_m,q_\alpha,q_\beta)$.
Now, recall that $Y$ is mapped isomorphically to $Z\subset \p^3$.
We obtain 
\dia
H^0(\n_{X/P_m}|_Y) & H^0(\os_Z(2m)) \ar[l,"\cong"] &  H^0(\on{3}(2m)) \ar[l,twoheadrightarrow] & H^0(\id_Z(2m)) \ar[l,hookrightarrow] \\
H^0(\n_{Y/P_m}) \ar[u,"\rho"]  & H^0(\os_{Z}(m)\oplus \os_Z(\alpha)\oplus\os_Z(\beta)) \ar[l,"\cong"] \ar[u] & H^0(\on{3}(m)\oplus \on{3}(\alpha)\oplus\on{3}(\beta)) \ar[l,twoheadrightarrow] \ar[u,"\rho'"] \\
\mma
where we used that the restriction maps from line bundles on $\p^3$ to $Z$ are surjective on global sections, see \autoref{numericsci}.
By commutativity, to prove that $\rho$ is surjective it is enough to show that $\rho'$ is such.
By our construction, the ideal of $Z$ is generated by the two polynomials $p_\alpha,p_\beta$; the image of the rightmost horizontal map is exactly the degree $2m$ part of this ideal.
Since $\pi$ is identified to the scalar product with $(-p_m,q_\alpha,q_\beta)$ we deduce that $\rho'$ is surjective if and only if the polynomials $p_\alpha,p_\beta,-p_m,q_\alpha,q_\beta$ generate all the homogeneous polynomials of degree $2m$.
But we have seen in \autoref{polinomi} that this holds if those five polynomials are general, equivalently if $(Y,X)$ is general in the variety $\mathscr{I}$ considered in \autoref{a=b}.

\textbf{Step 3: Conclusion}

\noindent The first formula for $h^0(\n_{Y/X})$ follows immediately from the surjectivity of $\rho$ and \eqref{normali}.
For the second one, we first use $h^0(\n_{Y/X})=\chi(\n_{Y/X})+h^1(\n_{Y/X})$ and then by Riemann-Roch we have
\[\chi(\n_{Y/X})=deg(\n_{Y/X})+2(1-p_a(Y))=deg(\omega_Z\otimes\omega_X\du)-deg(\omega_Z)=deg(\os_Z(4-m))=\alpha\cdot \beta\cdot (4-m)\]
by \autoref{divisorialcover} and the usual $Y\cong Z$.
For the computation of $h^1(\n_{Y/X})$, by \autoref{numericsci} we have $\omega_Y\cong \os_Y(\alpha+\beta-4)$.
It follows by Serre duality that $h^1(\os_Y(2m))=h^0(\omega_Y(-2m))=0$ being $\alpha,\beta\le m$; together with the surjectivity of $\rho$, it implies that $h^1(\n_{Y/X})=h^1(\n_{Y/P_m})$.
In a similar manner, noting that $Y\cong Z$ is connected, we have
\[h^1(\os_Y(\alpha))=\begin{cases} 0  \hfill \text{otherwise} \\
1 \qquad \hfill  \beta=m=4 \\
\end{cases}
\qquad h^1(\os_Y(\beta))=h^1(\os_Y(m))=\begin{cases} 0  \hfill \text{otherwise} \\
1 \qquad \hfill  \alpha=\beta=m=4 \\
\end{cases}\]
and the claim follows.
\end{proof}

\Co\label{smooth}
For $X,Y$ general, the point $[Y]$ is regular in the Hilbert scheme on $X$ if $\alpha\le \beta<4$.
\io
\begin{proof}
    From \autoref{cohonormale} we know that $\alpha\le \beta<4$ implies $h^1(\n_{Y/X})=0$ for general $Y\subset X$.
    Since we have seen that $Y\subset X$ is l.c.i. then, by \cite[Prop. 6.5.2]{FGAexplained} we conclude.
\end{proof}

\Le\label{HSrk2}
Let $(X,H)$ be a regular, $n$-dimensional polarised variety and $\e$ a rank $2$, initialised aCM bundle fitting 
\sesl{\os_X}{\e}{\id_Y(\gamma)}{solita'''}
If $n\ge3$ then 
\begin{itemize}
\item $hom(\e,\e)=h^0(X,\e\otimes\id_Y)+h^0(\e(-\gamma))$
    \item $ext^1(\e,\e)=h^{0}(X,\n_Y)-h^0(X,\e)+h^0(X,\e\otimes\id_Y)$
    \item $ext^j(\e,\e)=h^{j-1}(X,\n_Y)$ for all $1<j<n-1$.
\end{itemize} 
Moreover, if $K_X\sim kH$ with $k<0$, then $ext^{n-1}(\e,\e)=h^{n-2}(X,\n_Y)$  and $ext^{n}(\e,\e)=h^{n}(\e)$.
\ma
\begin{proof}
Recall that $\e|_Y\cong \n_Y$ hence, tensoring by $\e$ the standard sequence associated to $Y$ we get
   \ses{\e\otimes\id_Y}{\e}{\n_{Y/X}}
   Being $n\ge2$ and $\e$ aCM we can compute 
   \[h^{j}(X,\e\otimes\id_Y)=\begin{cases}
       h^{0}(\n_Y)-h^0(\e)+h^0(\e\otimes\id_Y) \qquad \hfill j=1\\
       h^{j-1}(\n_{Y/X}) \hfill 1<j<n\\
       h^n(\e) \hfill j=n. 
   \end{cases}\]
  Tensoring \eqref{solita'''} by $\e\du\cong \e(-\gamma)$ it reads
    \sesl{\e(-\gamma)}{\e\otimes\e\du}{\e\otimes\id_Y}{utd}
    Then, being $\e$ aCM we have $hom(\e,\e)=h^0(\e\otimes\id_Y)+h^0(\e(-\gamma))$ and being $\e$ locally free we get
    \[ext^j(\e,\e)=h^j(\e\otimes\e\du)= h^j(X,\e\otimes\id_Y)\]
    for $0<j<n-1$.
    In particular, since $n\ge3$ we obtain the formula for $ext^1(\e,\e)$.
    Moreover, if $K_X\sim kH$ with $k<0$ then, by Serre duality we have 
    \[h^n(\e(-\gamma))=h^n(\e\du)=h^0(\e(K_X))=0\]
    being $\e$ initialised, hence the above claim also holds for $j=n-1,n$.
\end{proof}

\Co\label{dimensioniacm}
Suppose that $char(\gk)\ne 2$ and $\gk$ infinite.
If $X$ is a general double covering of $\p^3$ ramified over a surface of degree $2m=4,6,8$.
Then, for any $1\le\alpha,\beta\le m$ there are rank $2$ aCM bundles $\e$ with $ext^2(\e,\e)=0$.
Moreover, assuming $\alpha+\beta>m$, they are stable and give regular points in a component of the expected dimension, equal to $ext^1(\e,\e)$, of their moduli spaces. 
\io
\begin{proof}
Recall that such $\e$ fits in the sequence \eqref{solita'''} where $Y$ is a curve as in \autoref{cohonormale}.
Moreover, from the uniqueness in \autoref{section-rk2} we deduce $h^0(\e\otimes\id_Y)=1$.
    For $2m=4,6$ the variety $X$ is Fano by \autoref{divisorialcover} so the vanishing of $ext^2(\e,\e)$ for a general such bundle follows directly from \autoref{HSrk2} and \autoref{cohonormale}. If $2m=8$ then $\omega_X\cong \os_X$ by \autoref{divisorialcover}.
 Using Serre duality and \autoref{HSrk2} we get
    \[ext^2(\e,\e)=ext^1(\e,\e)=h^{0}(X,\n_Y)-h^0(X,\e)+h^0(\e\otimes\id_Y)=h^{0}(X,\n_Y)-h^0(X,\e)+1.\]
    Recalling \autoref{cohonormale} we have $h^{0}(X,\n_Y)=h^{1}(X,\n_Y)$ since $m=4$.
    Finally, using \autoref{cohonormale} and \autoref{splittype} we conclude that $ext^2(\e,\e)=0$.
  
    The remaining part of the proof is uniform in $m$.
    Being $\alpha+\beta>m$ then \autoref{condstab} implies that $\e$ is slope stable, in particular those moduli spaces are well defined.
    Finally, \cite[Prop. 6.5.1]{FGAexplained} tells us that sheaves with $ext^2(\e,\e)=0$ give us regular points of their moduli spaces whose tangent space is $Ext^1(\e,\e)$.
\end{proof}
Note that for $\alpha=\beta=m=2$ there is a description of all the non-Ulrich sheaves appearing in this moduli space in \cite[Prop. A.1]{FLZ}.
The case $m=4$ is particularly interesting. 
\Co\label{sferici}
Suppose that $char(\gk)\ne 2$, $\gk$ infinite and pick any $1\le\alpha\le\beta\le4$ such that $\alpha+\beta>4$.
For any such choice of $\alpha,\beta$, on a general octic double solid $X$ have a spherical aCM bundle of rank $2$.
Moreover, if $\gk=\C$ then there are at least two non-isomorphic such bundles.
\io
\begin{proof}
From the above proof $ext^1(\e,\e)=0=ext^2(\e,\e)$ and $h^0(\e\otimes\id_Y)=1$.
Moreover, $\gamma=\alpha+\beta-4>0$ and $\e$ initialised, so that $h^0(\e(-\gamma))=0$.
Thus, by \autoref{HSrk2} $\e$ is simple and by Serre duality we have $1=hom(\e,\e)=ext^3(\e,\e)$, hence $\e$ is spherical.
Moreover, if $\gk=\C$ we have $Pic(X)\cong\Z$ by \autoref{Piccover} and for a general $X$ the Picard group of the ramification locus $R$ is cyclic hence $R$ cannot contain curves as the $Y$ above.
Therefore, by \autoref{fisso'}, $\e$ and $\iota^*\e$ are not isomorphic.
\end{proof}

Note that the Hilbert schemes of curves $Y$ and bundles $\e$ are linked by an actual morphism.
\Oss
The sequence \eqref{solita'''} is an instance of the Hartshorne--Serre correspondence between locally complete intersection codimension $2$ subschemes and rank $2$ vector bundles.
In case the resulting sheaves are semistable, in \cite[Cor. 2.12, 2.14]{tesi} it is shown that this correspondence glues to a morphism between the locus in the Hilbert scheme containing the curves $Y$ and the moduli space of sheaves containing $\e$.
Moreover, the fibre of this map over the point $[\e]$ parametrising a stable bundle is identified with $\p(H^0(\e))$, so this morphism is smooth on the stable locus.
In particular, for $\alpha\le\beta<4$ we could have computed the dimension of the moduli spaces containing $\e$ and deduce smoothness directly from \autoref{smooth}, using \cite[\href{https://stacks.math.columbia.edu/tag/02K5}{Tag 02K5}]{Stacks}.
On the other hand, suppose $\beta=4=m$.
Using the fact that composition of smooth morphism is smooth, we can deduce the smoothness of the Hilbert scheme of the curves $Y$ from smoothness of the moduli space of $\e$.
In particular, \autoref{smooth} holds even if $\beta=4$.
\one

\section{Applications to divisorial coverings of arbitrary varieties}
In this section, we first prove a variation of \cite[Thm. 1.3]{Casnati_wild}, see also \cite[Thm. 2.6]{EisSch} for the analogous statement regarding the (non fiber) product, which, together with the results of the previous part, allows us to construct Ulrich bundles on new classes of varieties. 
\Prop\label{pullback}
Suppose we have two equidimensional schemes $X_1,X_2$ and maps as in
\begin{equation}\label{pb}
\begin{tikzcd}
    X:=X_1\times_{\p^n}X_2 \ar[r, "p_2"] \ar[d, "p_1"] & X_2 \ar[d,twoheadrightarrow ,"f_2"] \\
X_1 \ar[r, "f_1"] & \p^n \\
\end{tikzcd}
\end{equation}
with $f_1,f_2$ finite and such that $f_2$ is also surjective.
If $\e_i$ is an Ulrich (aCM) sheaf on $(X_i,f_i^*\osn(1))$ then $p_1^*\e_1\otimes p_2^*\e_2$ is Ulrich (aCM) on $(X,p_1^*f_1^*\osn(1))$.
\One
\begin{proof}
We treat the Ulrich case, the aCM one is analogous.
    First of all, $p_1,p_2$ are pullbacks of finite morphisms and hence finite, so that $p_1\circ f_1=p_2\circ f_2$ is finite and $(X,p_1^*f_1^*\osn(1))$ is a polarised variety.
    Moreover, $p_1$ is surjective being $f_2$ such.
    By projection formula we have $(p_1)_*(p_1^*\e_1\otimes p_2^*\e_2)\cong \e_1\otimes (p_1)_*p_2^*(\e_2)$.
    Being $f_2$ finite, by base change \cite[Thm. 12.6]{GW1}, we have 
    \[(p_1)_*p_2^*(\e_2)\cong f_1^*(f_2)_*\e_2\cong f_1^*\osn^N\cong \os_{X_1}^N\]
    where we used \autoref{noether}.
    We conclude that $(p_1)_*(p_1^*\e_1\otimes p_2^*\e_2)\cong \e_1^N$ is Ulrich on $(X_1,f_1^*\osn(1))$ hence by \autoref{proiezione} we proved that $p_1^*\e_1\otimes p_2^*\e_2$ is Ulrich.  
\end{proof}

As an immediate corollary we have an easy proof of the following well-known fact.

\Co\label{restrizione}
Let $\e$ be an Ulrich (aCM) sheaf on an $n$-dimensional polarised scheme $(X_2,H_2)$ over an infinite field $\gk$.
Suppose $n\geq 1$.
For any $D\in |\os_{X_2}(1)|$ define $\e_D:=\e\otimes\os_D$ and $\os_D(1):=\os_{X_2}(1)|_D$, then $\e_D$ is Ulrich (aCM) on $(D,\os_D(1))$.
The same holds for any complete intersection of divisors in $|\os_{X_2}(1)|$.
\io
\begin{proof}
    In \autoref{pullback}, we can choose $X_1\cong \p^{n-1}\subset\p^n$ with the Ulrich sheaf $\os_{\p^{n-1}}$ such that $D=X:=X_1\times_{\p^n}X_2$.
    By arguing inductively, we conclude the second assertion. 
\end{proof}
Furthermore, in the aCM case, the same argument allows us to consider complete intersection of divisors in different $|\os_{X_2}(i)|$ but we need to ask that they lie in $f^*|\osn(i)|\subseteq |\os_{X_2}(i)|$.

A more interesting application is the following.
Consider $X_1\subset\p^n$ and set $\os_{X_1}(1):=\osn(1)|_{X_1}$.
In analogy with the case of $\p^n$, we define a finite covering $f:X\to X_1$ to be \textit{divisorial} if $f$ factors through a closed embedding $X\subset \p(\os_{X_1}\oplus\os_{X_1}(m))$, for some $m>0$, followed by the standard projection.
\Prop\label{divgeneral}
Let $(X_1,\os_{X_1}(1))$ be a polarised variety with an Ulrich sheaf $\e$.
Let $p_1:X\to X_1$ be a divisorial covering of degree $d$ and equation $\sum_{j=0}^dp_{d-j}t^j$ and define $\os_X(1):=p_1^*\os_{X_1}(1)$.
If there is some finite morphism $f_1:X_1\to\p^n$ such that $f_1^*\osn(1)\cong\os_{X_1}(1)$ and for all $j=0,\dots d$ we have $p_j\in f_1^*H^0(\osn(jm))\subseteq H^0(\os_{X_1}(jm))$ then $(X,\os_X(1))$ admits an Ulrich sheaf. 
\One
\begin{proof}
    Choose polynomials $\overline{p_j}\in H^0(\osn(jm))$ such that $p_j=f_1^*\overline{p_j}$.
    Call $f_2:X_2\to\p^n$ the divisorial cover of $\p^n$ defined by the polynomial $\sum_{j=0}^d\overline{p_{d-j}}t^j$.
    Then, the morphism $p_1:X\to X_1$ is the base change of $f_2$ along $f_1$ hence we can apply \autoref{pullback} and conclude. 
\end{proof}
\Co\label{pjn}
Any divisorial covering of a projectively normal variety admits Ulrich sheaves
\io

We conclude by applying this result to some special varieties.
First, some Fano $3$-folds where existence of Ulrich bundles is not yet covered by the existing theory.
\Co\label{i1}
The smooth Fano $3$-folds in the families \href{https://www.fanography.info/1-2}{1-2 b)} and \href{https://www.fanography.info/1-5}{1-5 b)} in \cite{fanography} admit Ulrich bundles.
\io
\begin{proof}
    The first family parametrises double covering of the smooth quadric $Q\subset\p^4$ branched along a divisor cut by a quartic hypersurface.
    Since $Q$ is projectively normal, we just apply \autoref{pjn}.
    
    The second family parametrises varieties constructed as follows.
    Call $G:=Gr(2,5)$.
    The linear system associated to the ample generator $\os_G(1)$ of $Pic(G)$ gives a projectively normal embedding in $\p^9$, called Pl\"ucker embedding, and $\omega_G\cong \os_G(-5)$.
    If we call $W$ a smooth intersection with a codimension $3$ subspace then $W$ is still projectively normal and $\on{9}(-2)|_W\cong \omega_W$.
    The family we search for is made of double coverings of such $W$-s branched along a divisor in $|\omega_X\du|$, in particular we can again apply \autoref{pjn}.
\end{proof}

Finally, we look at some Horikawa surfaces.
With the notations of $\autoref{pullback}$, let $f_1:X_1\hookrightarrow \p^3$ to be a regular quadric, $X_2\to\p^3$ to be a regular double covering branched along a degree $6$ surface $B$.
Then, $p_1:X\to X_1$ is a double covering branched along $f_1^*B$ and $X$ is as in \cite[Thm. 1.6 iii); case $d=0$, $n=3$]{Horikawa1}.
\Prop\label{horik}
Suppose $\gk$ infinite and $char(\gk)\ne2$.
Let $(X_1,\os_{X_1}(1))$ be a regular quadric surface with the standard polarisation and $D\subset X_1$ a general divisor in $|\os_{X_1}(6)|$. 
If $p_1:X\to X_1$ is the degree $2$ covering branched along $D$ then $(X,p_1^*\os_{X_1}(1))$  admits rank $2$ Ulrich sheaves.
\One
\begin{proof}
We know $Pic(X_1)\cong \Z^2$ and it is generated by the classes of two intersecting lines; it is easily verified that they give Ulrich line bundles on $(X_1,\os_{X_1}(1))$.
Note that a general $D$ as in the statement is cut on $X_1$ by a general divisor $D'\in |\on{3}(6)|$, since $X_1$ is projectively normal.
If $X_2$ is the double covering of $\p^3$ branched along $D'$, then it admits Ulrich sheaves of rank $2$ by \autoref{rk2esistenza}.
We are in the setting of \autoref{pullback}, hence we conclude.
\end{proof}

\newpage

\appendix
\section{}\label{appendix}

We work with polynomials in $4$ variable.
In the following, we choose five homogeneous polynomials of degrees $\alpha,2m-\alpha,\beta,2m-\beta,m$ and let \textit{Macaulay2} check that the ideal they generate contains all degree $2m$ polynomials.
I thank Fulvio Gesmundo for suggesting me this strategy.

If $char(\gk)=0$ then we first check the same statement for the ring $R:=\Q[x,y,z,w]$, then, since it is a surjectivity claim, it is enough to tensor with our field $\gk$ to get the thesis.
We will describe explicitly the steps for the case of Ulrich bundles on the octic double solid.
\begin{lstlisting}[language=Macaulay2]
--Define our polynomial ring
R=QQ[x,y,z,w]
--Consider the special ideal
I=ideal{x^4,y^4,z^4,w^4,(x+y+z+w)^4}
--We can compute a basis for the degree 8 part of the quotient ring R/I and check that 
--it is just the 0 vector 
basis(8,R/I)
 \end{lstlisting}

To get the curves $Y$ for all the values of $\alpha,\beta, m$ we consider the ideal 
\[I=I(\alpha,\beta,m):=(x^\alpha,y^\beta,z^{2m-\alpha},w^{2m-\beta},(x+y+z+w)^m)\]
and search for a basis of the degree $2m$ part of $R/I$.

Now lets suppose that $char(\gk)$ is arbitrary, but not $2$.
This time we should work with polynomials in $\Z[x,y,z,w]$ so that, if our claim is true then we can tensor by our favorite field and get the thesis.
It is enough to run the same procedure but with the term $(x+y+z+w)^m$ replaced by the following polynomial: expand $(x+y+z+w)^m$ as we were in characteristic $0$ and replace all its coefficients with a $1$, i.e. we get the sum of all degree $m$ monomials.